\documentclass[sn-mathphys,Numbered]{sn-jnl}

\usepackage{graphicx}%
\usepackage{multirow}%
\usepackage{amsmath,amssymb,amsfonts}%
\usepackage{amsthm}%
\usepackage{mathrsfs}%
\usepackage[title]{appendix}%
\usepackage{xcolor}%
\usepackage{textcomp}%
\usepackage{manyfoot}%
\usepackage{booktabs}%
\usepackage{algorithm}%
\usepackage{algorithmicx}%
\usepackage{algpseudocode}%
\usepackage{listings}%

\newcommand {\R}        {{\mathbb R}}

\newcommand {\Rnn}      {\R^{n \times n}}

\newcommand {\Rnm}      {\R^{n\times m}}

\newcommand {\Rmm}      {\R^{m\times m}}
\newcommand {\mat}      [1] {\left[\begin{array}{#1}}
\newcommand {\rix}          {\end{array}\right]}
\newcommand {\eproof}
      {\space
        {\ \vbox{\hrule\hbox{\vrule height1.3ex\hskip0.8ex\vrule}\hrule}}
        \par}
\newcommand {\eq} [1] {\begin{equation}\label{#1}}
\newcommand {\en} {\end{equation}}

\theoremstyle{thmstyleone}%
\newtheorem{theorem}{Theorem}
\theoremstyle{thmstyletwo}%
\newtheorem{example}{Example}%

\theoremstyle{thmstylethree}%

\begin{document}

\title[Article Title]{Exact deflation for  accurate SVD computation  of   nonnegative bidiagonal products of arbitrary rank}


\author[1]{\fnm{Huang} \sur{Rong}}\email{rongh98@hunnu.edu.cn}

\author[2]{\fnm{Jungong} \sur{Xue}}\email{xuej@fudan.edu.cn}

\equalcont{These authors contributed equally to this work.}

\affil[1]{\orgdiv{School of Mathematics and Statistics}, \orgname{Hunan Normal University}, \orgaddress{ \city{Changsha}, \postcode{410081}, \state{Hunan}, \country{China}}}

\affil[2]{\orgdiv{School of Mathematical Science}, \orgname{Fudan University}, \orgaddress{ \city{Shanghai}, \postcode{200433}, \state{Shanghai}, \country{China}}}



\abstract{ 


Dealing with zero singular values can be quite challenging, as they have the potential to cause numerous numerical difficulties. This paper presents a method for computing the singular value decomposition (SVD) of a nonnegative bidiagonal product of arbitrary rank, regardless of whether the factors are of full rank or rank-deficient, square or rectangular. A key feature of our method is its ability to exactly deflate all zero singular values with a favorable complexity, irrespective of rank deficiency and ill conditioning. Furthermore, it ensures the computation of nonzero singular values, no matter how small they may be, with high relative accuracy. Additionally, our method is well-suited for accurately computing the SVDs of arbitrary submatrices, leveraging an approach to extract their representations from the original product. We have conducted error analysis and numerical experiments to validate the claimed high relative accuracy.

 }




\maketitle

\section{Introduction}\label{sec1}

 Computing the   singular value decomposition (SVD) of matrix products or quotients is crucial in various applications, including statistical realization, control theory, canonical correlations, and source separation   \cite{Heath,How,ste}. Many algorithms have been developed to compute  SVDs for two matrices \cite{Bai1,Hua201, Moler}, three matrices \cite{Bo91,Chu1}, and even long matrix products \cite{MBar10,Chu,Golub,ste}. While these algorithms are backward stable and can compute SVDs with high absolute accuracy, they often struggle to  accurately compute tiny singular values, which are important in many applications \cite{Ar91,ste}, especially with high relative accuracy. Efforts have been made to design high-accuracy SVD algorithms for a single full-rank matrix \cite{Dem90,DemG99,Dem08,Dem99,dop09,dop03,Parlett,Hua191,Ye08}. However, these algorithms may not be directly applicable to   scenarios involving multiple matrices. Though not abundant, there are several algorithms designed   to accurately compute SVDs for two or three full-rank matrices,    as warranted by the matrix entries    \cite{DemG99,Drmac2}. However, the challenge remains   when it comes to computing SVDs  involving multiple matrices with high relative accuracy.  This is particularly true when dealing with rank deficiency, where the presence of zero singular values can cause numerous numerical difficulties.
 This paper specifically addresses the challenge of accurately computing the SVD of a nonnegative bidiagonal product of arbitrary rank, i.e.,
\eq{eq:biprod}A=B_1B_2\ldots B_K\in\mathbb{R}^{n_{0}\times n_{K}},\en
where  $B_k\in {\mathbb R}^{n_{k-1}\times n_{k}}, \ 1\le k\le K,$  are nonnegative lower or upper  bidiagonal  that are allowed to be rank-deficient and  rectangular. The proposed method aims to overcome the issues caused by rank deficiency and ill conditioning, ensuring that zero singular values are exactly identified and returned, while the nonzero singular values, regardless of their small magnitude, are computed with high relative accuracy.

As shown in Section 2 later, there are numerous scenarios where numerical computations involving the product form (\ref{eq:biprod}) are necessary. Indeed, such bidiagonal products of full rank are explicitly accommodated by various structured matrices, including  Cauchy, Vandermonde, Bernstein-Vandermonde and Cauchy-Vandermonde totally nonnegative (TN) matrices with distinct nodes \cite{kov05,Mar10,MarM10,Mar,JMar}. Utilizing the full-rank bidiagonal factors as input, one can apply the algorithms in the {\it TNTool} \cite{kov} to address a wide range of numerical linear algebra problems, including the SVD problem, with high relative accuracy. On the other hand, if the aforementioned structured matrices   contain repeated nodes, they   inherently become  rank-deficient.  Recently, an accurate algorithm has been provided for computing eigenvalues of a singular square TN matrix with high relative accuracy (returning zero eigenvalues exactly)  \cite{kov19}. However, as pointed out in \cite{kov19}, the accurate computation for the SVD problem of rank-deficient TN matrices  remains an open question. This paper addresses the challenging issue stemming from rank deficiency through our proposed method.


Dealing with zero singular values can indeed be a challenging issue in numerical computations, as they can lead to various numerical difficulties. One attractive aspect of our method is its ability to exactly deflate   zero singular values. Let us  provide a brief explanation of our approach. We utilize nontrivial element pairs within the factors $B_k$'s as  input. Interestingly, these elements inherently reveal the presence of zero singular values. For instance, consider the factor $B_{1} $  as in (\ref{eq:biprod}) that is lower bidiagonal  with its $(i,i-1)$th and $(i,i)$th elements being zero. This implies that the $i$th row of  $A$
 is zero, indicating a potential zero singular value that can be deflated by removing this row. The key challenge then lies in extracting the representation of the resulting submatrix   from the original one, a task which our proposed method adeptly addresses.

 It is worth noting that a methodology was initially presented for handling submatrix representations of full-rank TN matrices, which can be decomposed as a product of nonsingular bidiagonal factors \cite{kov07}. However, this approach encounters difficulties when confronted with matrices that exhibit rank deficiency.   A subsequent approach  addresses the computation of representations for matrices derived from   square TN matrices by removing a column and appending a zero column at the end \cite{kov19}. However, a gap exists in the ability to compute the representation  for arbitrary submatrices. Moreover, this capability is particularly vital when dealing with rank deficiency. Therefore, different from the aforementioned methods, we introduce the  updating and downdating procedures specifically tailored to extract the representation of arbitrary submatrices from the original one for any nonnegative bidiagonal product, including general TN matrices. The   extraction method offers a comprehensive solution that extends the capabilities of existing methods and enables the computation of submatrix representations in a broader range of general TN matrices.

 Basing on the extracting method, our  deflation method is designed to exactly deflate zero singular values. Let us briefly outline the deflation method. Given the   product $A=B_{1}B_{2}\ldots B_{K}\in {\mathbb R}^{n_{0}\times n_{K}}$ as in (\ref{eq:biprod}) having $S$ nontrivial element pairs in the factors $B_{k}\in {\mathbb R}^{n_{k-1}\times n_{k}}$ ($1\leq k\leq K$), one key idea   is to first split  the product   $A$ into two parts
 \[A=A_{2}A_{1},{\rm where}~ A_{1}=B_{T+1}\ldots B_{K}\in {\mathbb R}^{r\times n_{K}},~A_{2}=B_{1}\ldots B_{T}\in {\mathbb R}^{n_{0}\times r},\] by determining the minimum dimension among the factors, denoted as $r=\min_{0\leq k\leq K}\{n_{k}\}=n_{T}$. We  then compute the representations of $A_{1}$ and $A_{2},$  on which we perform the  deflation for $A=A_{2}A_{1}$  in the following process. This deflation process is based on the minimum dimension among the factors, ensuring that the deflation is carried out effectively. In what follows, we denote by
      $A(\alpha|:]$ (or $A[:|\alpha)$)   the submatrix of $A$   deleting  the rows (or columns) indexed in the integer array $\alpha$.

 \begin{itemize}\item 
First, deleting the zero rows indexed in $\alpha$ of $A_{1}$, which are revealed by its   representation, and removing
 the corresponding    columns of $A_{2},$ we get that
 $$
 A=A'_{2}A'_{1}:=A_{2}[:|\alpha)\cdot A_{1}(\alpha|:],
 $$
where  the   representations of $A'_{1}$ and $A'_{2}$   are extracted from those of $A_{1}$ and $A_{2}$, respectively. Then,    a nonsingular   bidiagonal matrix $L$ is constructed such that  zero rows of $A_{1}^{(1)}=LA'_{1}$ are further revealed  by its  representation, reformulating $A$ as
\[A=A''_{2}A_{1}^{(1)}:=(A_{2}'L^{-1})(LA'_{1}), \]
where the  representations of $A_1^{(1)}$ and $A''_2$ are extracted  from those of $A'_{1}$ and $A'_{2}$, respectively.
\item Second,  we perform an orthogonal transformation to deflate some zero singular values of $A$ by deleting  zero rows of $A''_{2}$ revealed by its  representation, arriving at
\[A^{(1)}=A^{(1)}_{2}A^{(1)}_{1}:=(G^{T}A''_{2})A_{1}^{(1)},\]
where $G$ are orthogonal, and the   representation of $A^{(1)}_{2}$ is extracted from that of $A''_{2}$.
\item Further,
via the resulting representation of $A^{(1)}_{2}$, we  delete the zero columns indexed in $\beta$ of $A^{(1)}_{2}$ and remove the corresponding   rows of $A^{(1)}_{1}$ to get that
\[A^{(1)}=A'^{(1)}_{2}A'^{(1)}_{1}:=A^{(1)}_{2}[:|\beta)\cdot A^{(1)}_{1}(\beta|:],\]
where the   representations of $A'^{(1)}_{1}$ and $A'^{(1)}_{2}$ are extracted from those of $A^{(1)}_{1}$ and $A^{(1)}_{2}$, respectively.
Then,     a nonsingular  bidiagonal matrix $U$ is constructed such that   zero columns of $A^{(2)}_{2}=A'^{(1)}_{2}U$ are further revealed  by its representation, reformulating $A^{(1)}$ as
\[A^{(1)}=A^{(2)}_{2}A''^{(1)}_{1}:=(A'^{(1)}_{2}U)(U^{-1}A'^{(1)}_{1}), \]
where the   representations of  $A''^{(1)}_{1}$ and $A^{(2)}_{2}$ are extracted from those of $A'^{(1)}_{1}$ and $A'^{(1)}_{2}$, respectively.
\item Finally,   
we perform an orthogonal transformation to deflate some zero singular values of $A^{(1)}$ by  deleting  zero columns of $A''^{(1)}_{1}$ revealed by its representation, arriving at
\[A^{(2)}=A^{(2)}_{2}A^{(2)}_{1}:=A^{(2)}_{2}(A''^{(1)}_{1}V),\]
where $V$ is orthogonal, and the   representation of $A^{(2)}_{1}$ is extracted from that of $A''^{(1)}_{1}$.
\end{itemize}
We subsequently  perform the detection and deflation operations described above on the representations of  $A^{(2)}_{2}$ and $A^{(2)}_{1}$ to yeild  a product of a nonsingular diagonal matrix and a bidiagonal matrix of full rank.  Thereby, all the zero singular values  are deflated. Subsequently, the nonzero singular values  are computed by applying the available SVD algorithm   \cite{Dem90,Parlett,Parlett1}   to   the resulting full-rank bidiagonal matrix. The whole process  involves no subtraction of same-signed numbers, ensuring that zero singular values are exactly deflated and the  nonzero one are computed with high relative accuracy. Moreover, the overall computational cost is  at most $O(rS+\max\{n^{2}_{0}r,n^{2}_{K}r\})$, which  is preferable when $r\ll \min\{n_{0},n_{K}\}$. 

Remark that  a periodic  reduction method has been developed to   reduce a full-rank product of  structured matrices into a tridiagonal form, facilitating the accurate eigenvalue computation \cite{Hua18}. In addition, a reduction method has been provided to transform a full-rank product of an  upper triangular matrix and a  BF matrix into a bidiagonal form, tailored for the accurate generalized singular value computation  \cite{Hua191}. Notably, our  method here focuses on exactly deflating zero singular values of a rank-deficient bidiagonal product.

   The rest of the paper is organized as follows. In Section 2, we  illustrate the crucial role played by the bidiagonal product form when dealing with  rank deficiency and submatrices.   Section 3 introduces  a new method    to  extract the representation of an arbitrary submatrix from the original one, which is indispensable when  deflating zero singular values.    In Section 4,  we establish a   method     to  exactly deflate   zero singular values of    a rank-deficient nonnegative bidiagonal product in a preferable complexity. The overall SVD algorithm is designed to accurately compute  the SVD of such a product as well as that of its arbitrary submatrix.     Numerical experiments are performed in Section 5 to confirm the claimed high relative accuracy.

   Before proceeding, we present some notation and symbols.  Let $\alpha$ and $\beta$ be integer arrays.
      $A[\alpha|\beta]$ (or $A(\alpha|\beta)$) represents the submatrix of $A\in\Rnm$ having (or deleting) the rows and columns indexed in $\alpha$ and $\beta,$ respectively, and it is abbreviated to $A[\alpha]$ (or $A(\alpha)$) if $\alpha=\beta.$    The $n\times m$ identity matrix is denoted by $I_{n,m},$ and $I_{n,n}$ is abbreviated to $I_{n}$. In the sequel,  we mean by
 \[L_{f:s}={\bf bilow}(\{\bar{x}_{i},x_{i}\}_{i=f}^{s})\in\Rnn, ~U_{f:s}={\bf biupp}(\{\bar{y}_{i},y_{i}\}_{i=f}^{s})\in\Rnn \]
 the lower (upper) bidiagonal matrix obtained from  $I_{n}$ by replacing its $(i,i)$th and $(i+1,i)$th ($(i,i)$th and $(i,i+1)$th) entries by $\bar{x}_{i}$ and $x_{i}$ ($\bar{y}_{i}$ and $y_{i}$) for all $f\leq i\leq s$, respectively; with the convention $\{x_{n}\}=\emptyset$ ($\{y_{n}\}=\emptyset$) if $s=n$. All $\{\bar{x}_{i},x_{i}\}$ ($\{\bar{y}_{i},y_{i}\}$) are  called nontrivial {\it element pairs} of $L_{f:s}$ ($U_{f:s}$).   Set $L_{f:s}=I_{n}$ ($U_{f:s}=I_{n}$) if $f>s$ or $f>n$. Observe that $L_{f:s}=L_{f:f}L_{f+1:f+1}\ldots L_{s:s}$ and $U_{f:s}=U_{s:s}\ldots U_{f+1:f+1} U_{f:f}$. 

\section{Bidiagonal Product Forms} \label{sec2} In this section, we  underscore the crucial role played by the bidiagonal product form when dealing with the SVD computations for rank deficiency and submatrices.

\subsection{Computing SVDs of rank-deficient structured  matrices} \label{sec2.1}Many structured matrices with distinct nodes, such as Cauchy, Vandermonde, Cauchy-Vandermonde  and Bernstein-Vandermonde matrices, are full-rank and have been shown to have the explicit  bidiagonal decomposition   \cite{kov05,Mar10,MarM10,Mar,JMar}:
\eq{bd}A=L_{n-1}\ldots   L_{1}DU_{1}\ldots U_{m-1}\in\Rnm,\en 
 where $D\in {\mathbb R}^{n\times m}$ is diagonal, $L_{k}\in {\mathbb R}^{n\times n}$ (or $U_{k}\in {\mathbb R}^{m\times m}$) is   nonsingular lower (or upper)  bidiagonal. Based on the   bidiagonal decomposition, accurate SVD computations   have been performed via the available algorithms \cite{kov07}. However, the presence of repeated nodes in these structured matrices often leads to rank deficiency, presenting a significant obstacle to accurately computing their bidiagonal decompositions \cite{kov19}.   Recent efforts have  addressed this issue by exploiting special structures of the factors $L_k$
 in the bidiagonal decomposition  in the context of Vandermonde and    ($q$-, $h$-) Bernstein-Vandermonde matrices \cite{{Delg1}}. However, a general approach for representing  structured matrices with repeated nodes in a bidiagonal product form has remained elusive. In response to this gap, we  propose a systematic methodology that allows for the representation of these matrices with repeated nodes in a more manageable bidiagonal product form.

 In fact, these structured matrices containing repeated nodes can be effectively expressed as a product of bidiagonal matrices and structured matrices of the same type but with distinct nodes. This process is demonstrated as follows. 
Let the $(i,j)$-th entry of the structured matrix $A\in\Rnm$  be a function of $x_i$ and $y_j.$ The nodes $x_i,$ $1\le i\le n$ and $y_j,$ $1\le j\le m,$ are arranged such that the same nodes appear consecutively. More specifically, we let
\[\begin{cases} x_{i_{k}}=x_{i_{k}+1}=\ldots=x_{i_{k+1}-1}\neq x_{i_{k+1}},~\forall~1\leq k\leq n_{1},\cr y_{j_{k}}=y_{j_{k}+1}=\ldots=y_{j_{k+1}-1}\neq y_{j_{k+1}},~\forall~1\leq k\leq m_{1}, \end{cases}\] with the conventions $i_{1}=j_{1}=1$, $i_{n_{1}+1}=n+1$ and $j_{m_{1}+1}=m+1.$ Then the corresponding rows (and columns) of $A$ are repeated, i.e.,
\[\begin{cases}A[i_{k}|:]=A[i_{k}+1|:]=\ldots=A[i_{k+1}-1|:]\neq A[i_{k+1}|:],~\forall~1\leq k\leq  n_{1},\cr A[:|j_{k}]=A[:|j_{k}+1]=\ldots=A[:|j_{k+1}-1]\neq A[:|j_{k+1}],~ \forall~1\leq k\leq m_{1}.\end{cases} \]  A key observation is that $A$ can be obtained by expanding its submatrix $$A[i_{1},\ldots,i_{n_{1}}|j_{1},\ldots,j_{m_{1}}],$$ which is the structured matrix of the same type  associated with  distinct nodes,  via multiplication by rectangular identity matrices and the   bidiagonal matrices $E_{i}(x,y).$ Here, $E_{i}(x,y)$ is the matrix obtained from the identity $I_{n}$ by replacing its $(i,i)$ and $(i+1,i)$ entries with $x$ and $y$, respectively. We have
\begin{eqnarray*}A&=&\left(\prod_{k=1}^{n_{1}}\prod^{t=i_{k}}_{i_{k+1}-2}E_{t}(1,1)\right)
	\left(\prod_{k=2}^{n_{1}}\prod^{t=k}_{i_{k}-1}E_{t}(0,1)\right) \\&&I_{n,n_{1}}\cdot A[i_{1},\ldots,i_{n_{1}}|j_{1},\ldots,j_{m_{1}}]\cdot I_{m_{1},m}\nonumber \\ && \left(\prod_{k=2}^{m_{1}}\prod^{t=k}_{j_{k}-1}E_{t}(0,1)\right)^{T}\left(\prod_{k=1}^{m_{1}}\prod^{t=j_{k}}_{j_{k+1}-2}E_{t}(1,1)\right)^{T}.\nonumber \end{eqnarray*}
As long as the full-rank  bidiagonal decomposition $$
A[i_{1},\ldots,i_{n_{1}}|j_{1},\ldots,j_{m_{1}}]=L_{n_{1}-1}\ldots L_{1}DU_{1}\ldots U_{m_{1}-1}
$$ is available,
    $A$ can be expressed in the {\it bidiagonal product} form
\begin{eqnarray}A&=&\left(\prod_{k=1}^{n_{1}}\prod^{t=i_{k}}_{i_{k+1}-2}E_{t}(1,1)\right)
	\left(\prod_{k=2}^{n_{1}}\prod^{t=k}_{i_{k}-1}E_{t}(0,1)\right)\label{bfs}\\&&I_{n,n_{1}}\cdot L_{n_{1}-1}\ldots L_{1}DU_{1}\ldots U_{m_{1}-1}\cdot I_{m_{1},m}\nonumber \\ && \left(\prod_{k=2}^{m_{1}}\prod^{t=k}_{j_{k}-1}E_{t}(0,1)\right)^{T}\left(\prod_{k=1}^{m_{1}}\prod^{t=j_{k}}_{j_{k+1}-2}E_{t}(1,1)\right)^{T}.\nonumber \end{eqnarray}
This approach offers a unified framework for accurately computing the SVDs  involving   structured matrices with repeated nodes, as demonstrated  in Section \ref{sec5}. 

\subsection{Computing SVDs of submatrices} \label{sec2.2}The SVD computation of submatrices has garnered substantial interest from various applications  \cite{Hen19}. 
We interestingly  find out that an  arbitrary submatrix can be expressed as a product of  rectangular bidiagonal matrices and the original matrix. Let us concretely demonstrate this fact as follows.

Consider the submatrices $A(r|:]$ and $A[:|r)$ of $A\in\Rnm$. Denote
\[\mathcal{U}_{r:n-1}={\bf biupp}(\{0,1\}_{i=r}^{n-1})\in\mathbb{R}^{(n-1)\times n},\quad \mathcal{L}_{r:m-1}={\bf bilow}(\{0,1\}_{i=r}^{m-1})\in\mathbb{R}^{m\times (m-1)}.\] Then
\[A(r|:]=\mathcal{U}_{r:n-1}A,\quad A[:|r)=A\mathcal{L}_{r:m-1}.\]
Thereby, for an arbitrary submatrix $A(\alpha|\beta)$ of $A$, it can be expressed as
\[A(\alpha|\beta)= \prod_{i=1}\mathcal{U}_{\alpha_{i}:n-k+i-1} \cdot A\cdot \prod_{l}^{i=1}\mathcal{L}_{\beta_{i}:m-k+i-1}.\]
Here,
$\alpha=(\alpha_{i})\in Q_{k,n}$ and $\beta=(\beta_{i})\in Q_{l,m}$, and $Q_{k,n}$ denotes the set of the vectors of length $k$ consisting of distinctive integers from $\{1,\ldots,n\}.$
 In particular, if the original matrix $A$ has a bidiagonal decomposition, e.g., the form (\ref{bd})  is available,  then the submatrix $A(\alpha|\beta)$ is of the following bidiagonal product form
\[A(\alpha|\beta)= \prod_{i=1}\mathcal{U}_{\alpha_{i}:n-k+i-1} \cdot \prod_{n-1}^{i=1}L_{i}\cdot D\cdot \prod_{i=1}^{m-1}U_{i}\cdot \prod_{l}^{i=1}\mathcal{L}_{\beta_{i}:m-k+i-1},\]
As shown in Section \ref{sec5},  accurate SVD computations are achieved for arbitrary submatrices of many structured matrices as well as their products, by virtue of depicting submatrices as bidiagonal product forms.

Therefore, our primary focus in this paper is to conduct the SVD computation  for a bidiagonal product of arbitrary rank. Our proposed method aims to   compute with high relative accuracy the SVD of this product.



\section{The extracting method for arbitrary submatrices} In this section, we present a  method for extracting the representation of an arbitrary submatrix  from the original one in the context of a nonnegative bidiagonal product, specifically designed to handle rank deficiency by deflating zero singular values later.

Before proceeding, we introduce some notations that are frequently utilized later. By extending the representation  of a square TN matrix in \cite{kov19} to a rectangular   bidiagonal product, throughout the paper we always mean by
\eq{rp}A=:(\{\bar{g}_{ij},g_{ij}\})\in\Rnm\en
the    matrix   generated from the $n\times m$  element pairs $\{\bar{g}_{ij},g_{ij}\}_{i,j=1}^{n,m}$, which are stored as the matrix  $\mathcal{BD}(A)=(\{\bar{g}_{ij},g_{ij}\})\in\Rnm$, by the following bidiagonal product form \eq{gbf}A=L_{n-1}\ldots   L_{1}DU_{1}\ldots U_{m-1},\en
 where $D={\bf diag}(\{g_{ii}\}_{i=1}^{\min\{n,m\}})\in\Rnm$,  and
 \[\begin{cases}L_{k}={\bf bilow}(\{\bar{g}_{i+1,i+1-k},g_{i+1,i+1-k}\}_{i=k}^{\min\{n-1,m+k-1\}})\in \Rnn,~\forall~1\leq k\leq n-1,\cr
 U_{l}={\bf biupp}(\{\bar{g}_{i+1-l,i+1},g_{i+1-l,i+1}\}_{i=l}^{\min\{m-1,n+l-1\}})\in\Rmm,~\forall~1\leq l\leq m-1.\end{cases}\]
 We call   $\mathcal{BD}(A)$ to be the representation of $A$. Notice that  all $\bar{g}_{ii}$ has no impact on forming $A.$ In addition, by the  form (\ref{gbf}), without loss of generality, we always assume that $g_{ii}\neq 0$ for all $i\neq \min\{n,m\}$. If all $\bar{g}_{ij}\geq0$ and $g_{ij}\geq 0$, we denote $\mathcal{BD}(A)\geq 0$. Clearly, the representation of   $A^{T}$ is the transpose of $\mathcal{BD}(A)$, i.e.,
\eq{c1.2}A^{T}=:(\{\bar{g}_{ji},g_{ji}\})\in \mathbb{R}^{m\times n}.\en
Later, we will show that any nonnegative  bidiagonal product has the representation.

Given  the representation $\mathcal{BD}(A)=( \{\bar{g}_{ij},g_{ij}\})\in\Rnm$, we  provide the updating and downdating procedures  when appending and deleting  rows  (columns) of $A$.

 \begin{itemize}\item Suppose  $A'\in \mathbb{R}^{t\times m}$ ($t>n$) is obtained from $A\in\Rnm$ by attaching   $(t-n)$  zero rows to the last rows of $A.$  Then \[A'=I_{t,n}A=\mat{cc}A\\0_{(t-n)\times m}\rix= \mat{cc}L_{n-1}&\\&I_{t-n}\rix \ldots   \mat{cc}L_{1}&\\&I_{t-n}\rix \mat{cc}D\\0_{(t-n)\times m}\rix U_{1}\ldots U_{m-1}.\]
 Thus, $\mathcal{BD}(A')=(\{\bar{g}'_{ij},g'_{ij}\})\in \mathbb{R}^{t\times m}$ is obtained as
 \eq{dr2}\{\bar{g}'_{ij},g'_{ij}\}=\begin{cases}  (1,0),~\forall~n+1\leq i\leq t,~1\leq j\leq m,\cr \{\bar{g}_{ij},g_{ij}\},~{\rm otherwise}.\end{cases}\en
 \item Suppose $A'\in \mathbb{R}^{t\times m}$  ($t<n$) is obtained from $A\in\Rnm$ by deleting the last $(n-t)$ rows. Then \eq{drp}A'=I_{t,n}A=A[1:t|:]= L_{n-1}[1:t]\ldots   L_{1}[1:t] D[1:t|:]U_{1}\ldots U_{m-1}.\en Thus, $\mathcal{BD}(A')=(\{\bar{g}'_{ij},g'_{ij}\})\in \mathbb{R}^{t\times m}$ is obtained as
\eq{dr}\{\bar{g}'_{ij},g'_{ij}\}=\begin{cases}  \{\bar{g}_{tj},g_{tj} \prod_{k=1}^{j}\bar{g}_{t+1,k}\},~\forall~i=t,~1\leq j\leq \min\{t,m\},\cr \{\bar{g}_{ij},g_{ij}\},~{\rm otherwise},\end{cases}\en   costing  $2\cdot \min\{t,m\}$ subtraction-free  operations.  We illustrate the downdating procedure (\ref{dr}) by treating  the cases $t>m$ and $t\leq m$ separately with the aid of the  example $\mathcal{BD}(A)=(\{\bar{g}_{ij},g_{ij}\})\in \mathbb{R}^{7\times 4}$ and $t=5,2$ as follows.
\begin{itemize} \item For $t=5>m=4$, we have by (\ref{drp}) that
\begin{eqnarray*}A'&=&\tiny{\begin{bmatrix}\begin{smallmatrix}1&&&&\\&1&&&\\&&1&&\\&&&1&\\&&&&\bar{g}_{61}
\end{smallmatrix}\end{bmatrix}\begin{bmatrix}\begin{smallmatrix}1&&&&\\&1&&&\\&&1&&\\&&&\bar{g}_{51}&\\&&&g_{51}&\bar{g}_{62}
\end{smallmatrix}\end{bmatrix}\begin{bmatrix}\begin{smallmatrix}1&&&&\\&1&&&\\&&\bar{g}_{41}&&\\&&g_{41}&\bar{g}_{52}&\\&&&g_{52}&\bar{g}_{63}
\end{smallmatrix}\end{bmatrix}\begin{bmatrix}\begin{smallmatrix}1&&&&\\&\bar{g}_{31}&&&\\&g_{31}&\bar{g}_{42}&&\\&&g_{42}&\bar{g}_{53}&\\&&&g_{53}&\bar{g}_{64}
\end{smallmatrix}\end{bmatrix}\begin{bmatrix}\begin{smallmatrix}\bar{g}_{21}&&&&\\g_{21}&\bar{g}_{32}&&&\\&g_{32}&\bar{g}_{43}&&\\&&g_{43}&\bar{g}_{54}&\\&&&g_{54}&\bar{g}_{65}
\end{smallmatrix}\end{bmatrix}}\\
&&\tiny{\begin{bmatrix}\begin{smallmatrix}g_{11}&0&0&0\\0&g_{22}&0&0\\0&0&g_{33}&0\\0&0&0&g_{44}\\0&0&0&0&
\end{smallmatrix}\end{bmatrix}\begin{bmatrix}\begin{smallmatrix}\bar{g}_{12}&g_{12}&&\\&\bar{g}_{23}&g_{23}&\\&&\bar{g}_{34}&g_{34}\\&&&1
\end{smallmatrix}\end{bmatrix}\begin{bmatrix}\begin{smallmatrix}1&&&\\&\bar{g}_{13}&g_{13}&\\&&\bar{g}_{24}&g_{24}\\&&&1
\end{smallmatrix}\end{bmatrix}\begin{bmatrix}\begin{smallmatrix}1&&&\\&1& &\\&&\bar{g}_{14}&g_{14}\\&&&1
\end{smallmatrix}\end{bmatrix}}\\
&=&\tiny{\begin{bmatrix}\begin{smallmatrix}1&&&&\\&1&&&\\&&1&&\\&&&\bar{g}_{51}&\\&&&g_{51}\bar{g}_{61}&1
\end{smallmatrix}\end{bmatrix}\begin{bmatrix}\begin{smallmatrix}1&&&&\\&1&&&\\&&\bar{g}_{41}&&\\&&g_{41}&\bar{g}_{52}&\\&&&g_{52}\prod\limits_{k=1}^{2}\bar{g}_{6k}&1
\end{smallmatrix}\end{bmatrix}\begin{bmatrix}\begin{smallmatrix}1&&&&\\&\bar{g}_{31}&&&\\&g_{31}&\bar{g}_{42}&&\\&&g_{42}&\bar{g}_{53}&\\&&&g_{53}\prod\limits_{k=1}^{3}\bar{g}_{6k}&1
\end{smallmatrix}\end{bmatrix}\begin{bmatrix}\begin{smallmatrix}\bar{g}_{21}&&&&\\g_{21}&\bar{g}_{32}&&&\\&g_{32}&\bar{g}_{43}&&\\&&g_{43}&\bar{g}_{54}&
\\&&&g_{54}\prod\limits_{k=1}^{4}\bar{g}_{6k}&1
\end{smallmatrix}\end{bmatrix}}\\
&&\tiny{\begin{bmatrix}\begin{smallmatrix}g_{11}&0&0&0\\0&g_{22}&0&0\\0&0&g_{33}&0\\0&0&0&g_{44}\\0&0&0&0&
\end{smallmatrix}\end{bmatrix}\begin{bmatrix}\begin{smallmatrix}\bar{g}_{12}&g_{12}&&\\&\bar{g}_{23}&g_{23}&\\&&\bar{g}_{34}&g_{34}\\&&&1
\end{smallmatrix}\end{bmatrix}\begin{bmatrix}\begin{smallmatrix}1&&&\\&\bar{g}_{13}&g_{13}&\\&&\bar{g}_{24}&g_{24}\\&&&1
\end{smallmatrix}\end{bmatrix}\begin{bmatrix}\begin{smallmatrix}1&&&\\&1& &\\&&\bar{g}_{14}&g_{14}\\&&&1
\end{smallmatrix}\end{bmatrix}}\\
&&=:\tiny{\begin{bmatrix}
\begin{smallmatrix}\{\bar{g}_{11},g_{11}\}&\{\bar{g}_{12},g_{12}\}&\{\bar{g}_{13},g_{13}\}&\{\bar{g}_{14},g_{14}\}\\
\{\bar{g}_{21},g_{21}\}&\{\bar{g}_{22},g_{22}\}&\{\bar{g}_{23},g_{23}\}&\{\bar{g}_{24},g_{24}\}\\
\{\bar{g}_{31},g_{31}\}&\{\bar{g}_{32},g_{32}\}&\{\bar{g}_{33},g_{33}\}&\{\bar{g}_{34},g_{34}\}\\
\{\bar{g}_{41},g_{41}\}&\{\bar{g}_{42},g_{42}\}&\{\bar{g}_{43},g_{43}\}&\{\bar{g}_{44},g_{44}\}\\
\{\bar{g}_{51},g_{51}\bar{g}_{61}\}&\{\bar{g}_{52},g_{52}\prod\limits_{k=1}^{2}\bar{g}_{6k}\}&\{\bar{g}_{53},g_{53}
\prod\limits_{k=1}^{3}\bar{g}_{6k}\}&\{\bar{g}_{54},g_{54}\prod\limits_{k=1}^{4}\bar{g}_{6k}\}
\end{smallmatrix}\end{bmatrix}}.\end{eqnarray*}
\item For $t=2<m=4$, we have by (\ref{drp}) that
\begin{eqnarray*}A'&=&\begin{bmatrix}\begin{smallmatrix}1&\\&\bar{g}_{31}
\end{smallmatrix}\end{bmatrix}\begin{bmatrix}\begin{smallmatrix} \bar{g}_{21}&\\ g_{21}&\bar{g}_{32}
\end{smallmatrix}\end{bmatrix}\begin{bmatrix}\begin{smallmatrix}g_{11}&0&0&0\\0&g_{22}&0&0\end{smallmatrix}\end{bmatrix}\begin{bmatrix}\begin{smallmatrix}\bar{g}_{12}&g_{12}&&\\&\bar{g}_{23}&g_{23}&\\&&\bar{g}_{34}&g_{34}\\&&&1
\end{smallmatrix}\end{bmatrix}\begin{bmatrix}\begin{smallmatrix}1&&&\\&\bar{g}_{13}&g_{13}&\\&&\bar{g}_{24}&g_{24}\\&&&1
\end{smallmatrix}\end{bmatrix}\begin{bmatrix}\begin{smallmatrix}1&&&\\&1& &\\&&\bar{g}_{14}&g_{14}\\&&&1
\end{smallmatrix}\end{bmatrix}\\
&=&\begin{bmatrix}\begin{smallmatrix} \bar{g}_{21}&\\ g_{21}\bar{g}_{31}&1
\end{smallmatrix}\end{bmatrix}\begin{bmatrix}\begin{smallmatrix}g_{11}&0&0&0\\0&g_{22}\prod\limits_{r=1}^{2}\bar{g}_{3r}&0&0\end{smallmatrix}\end{bmatrix}
\begin{bmatrix}\begin{smallmatrix}\bar{g}_{12}&g_{12}&&\\&\bar{g}_{23}&g_{23}&\\&&1&0\\&&&1
\end{smallmatrix}\end{bmatrix}\begin{bmatrix}\begin{smallmatrix}1&&&\\&\bar{g}_{13}&g_{13}&\\&&\bar{g}_{24}&g_{24}\\&&&1
\end{smallmatrix}\end{bmatrix}\begin{bmatrix}\begin{smallmatrix}1&&&\\&1& &\\&&\bar{g}_{14}&g_{14}\\&&&1
\end{smallmatrix}\end{bmatrix}\\
&&=:\begin{bmatrix}
\begin{smallmatrix}\{\bar{g}_{11},g_{11}\}&\{\bar{g}_{12},g_{12}\}&\{\bar{g}_{13},g_{13}\}&\{\bar{g}_{14},g_{14}\}\\
\{\bar{g}_{21},g_{21}\bar{g}_{31}\}&\{\bar{g}_{22},g_{22}\prod\limits_{k=1}^{2}\bar{g}_{3k}\}&\{\bar{g}_{23},g_{23}\}&\{\bar{g}_{24},g_{24}\}
\end{smallmatrix}\end{bmatrix}.\end{eqnarray*}
\end{itemize}
 \end{itemize}
For the case of appending (or deleting) the columns of $A\in\Rnm$, i.e., $A'=AI_{m\times t}$ for $t>m$ and $t<m$,  the representation of $A'$ is obtained by applying (\ref{dr2}) (or (\ref{dr})) to  $I_{t\times m}A^{T}$ with the transposing operation (\ref{c1.2}).

In addition,     the command {\tt STNAddToPrevious} ({\tt STNAddToNext})  has been provided in  \cite{kov19}  to compute the representation of   $AE_{i}(x,y)$ ($E_{i}(x,y)A$), where $A$ is a square TN matrix and  $E_{i}(x,y)$ is an elementary matrix. Now, given  a rectangular   representation $  \mathcal{BD}(A)\in\Rnm$ and an   upper (lower) bidiagonal matrix $U$ ($L$), we extend the commands to compute the   representation  of $UA$ ($LA$)   by using the following passing-through operations, whose proofs are given in the appendix section.
\begin{itemize}\item Given
\[U_{1:n}={\bf biupp}(\{\bar{y}_{i},y_{i}\}_{i=1}^{n})\geq 0,~ L_{1:n}={\bf bilow}(\{\bar{x}_{i},x_{i}\}_{i=1}^{n})\geq 0, \]
then \eq{pt1}U_{1:n}L_{1:n}=\bar{L}_{1:n}\bar{U}_{1:n},\en
where
\[\bar{U}_{1:n}={\bf biupp}(\{\bar{y}'_{i},y'_{i}\}_{i=1}^{n}),~ \bar{L}_{1:n}={\bf bilow}(\{\bar{x}'_{i},x'_{i}\}_{i=1}^{n}) \]
satisfies that
for $1\leq i\leq n$, \[  \begin{cases}z_{1}=\bar{y}_{1}\bar{x}_{1},~w_{i}=z_{i}+x_{i}y_{i},\cr
\begin{cases}\bar{x}'_{i}=1,~\bar{y}'_{i}=w_{i},~x'_{i}=\frac{\bar{y}_{i+1}x_{i}}{\bar{y}'_{i}},~y'_{i}=\frac{y_{i}\bar{x}_{i+1}}{\bar{x}'_{i}},
~z_{i+1}=\frac{\bar{y}_{i+1}\bar{x}_{i+1}z_{i}}{\bar{x}'_{i}\bar{y}'_{i}},~{\rm if}~w_{i}\neq 0,\cr
\bar{x}'_{i}=0,~\bar{y}'_{i}= 1,~x'_{i}=\bar{y}_{i+1}x_{i},~y'_{i}=0,~z_{i+1}=\bar{y}_{i+1}\bar{x}_{i+1},~{\rm if}~w_{i}= 0~{\rm and}~y_{i}=0,\cr
\bar{x}'_{i}=1,~\bar{y}'_{i}=0,~x'_{i}=0,~y'_{i}=y_{i}\bar{x}_{i+1},~z_{i+1}=\bar{y}_{i+1}\bar{x}_{i+1},~{\rm otherwise},\end{cases}  \end{cases}
\] which costs at most $8n$ subtraction-free arithmetic operations.
\item Given
\[D={\rm diag}(d_{i})\in\Rnm,~U_{1:n}={\bf biupp}(\{\bar{y}_{i},y_{i}\}_{i=1}^{n}),\] then \eq{pt2}U_{1:n}D=\bar{D}\bar{U}_{1:m},~{\rm where}~\begin{cases} \bar{D}={\bf diag}(\bar{d}_{i})\in\Rnm,\cr \bar{U}_{1:m}={\bf biupp}(\{\bar{y}'_{i},y'_{i}\}_{i=1}^{m})\in\Rmm,\end{cases} \en
here, by a straightforward calculation,
\[\begin{cases} \begin{cases}\bar{d}_{i}=d_{i}\bar{y}_{i},~\bar{y}'_{i}=1,~y'_{i}=\frac{d_{i+1}y_{i}}{\bar{d}_{i}},~{\rm if}~ d_{i}\bar{y}_{i}\neq 0,\cr
		\bar{d}_{i}=1,~\bar{y}'_{i}=0,~y'_{i}=d_{i+1}y_{i},~{\rm if}~ d_{i}\bar{y}_{i}= 0,\end{cases}~\forall~1\leq i\leq \min\{n,m\},\cr
	\bar{y}'_{i}=1,~y'_{i}=0,~\forall~\min\{n,m\}<i\leq m,\end{cases}  \]  which costs at most $3\min\{n,m\}$ subtraction-free arithmetic operations.
\item  Given
 \[  U_{1:m}={\bf biupp}(\{\bar{y}_{i},y_{i}\}_{i=1}^{m})\geq 0,~  U'_{1:m}={\bf biupp}(\{\bar{x}_{i},x_{i}\}_{i=1}^{m})\geq 0,\] then
 \eq{pt3}U_{1:m}U'_{1:m}=\bar{U}'_{1:m}\bar{U}_{2:m},\en
 where \[ \bar{U}'_{1:m}={\bf biupp}(\{\bar{x}'_{i},x'_{i}\}_{i=1}^{m}),~
 \bar{U}_{2:m}={\bf biupp}(\{\bar{y}'_{i},y'_{i}\}_{i=2}^{m}),\]
satisfies that
for all $1\leq i\leq m-1$, \[\begin{cases}\bar{x}'_{1}=\bar{x}_{1}\bar{y}_{1},~ z_{1}=\bar{y}_{1}x_{1},~w_{i}=z_{i}+\bar{x}_{i+1}y_{i},\cr\begin{cases}
 \bar{y}'_{i+1}=1,~x'_{i}=w_{i},~\bar{x}'_{i+1}=\frac{\bar{x}_{i+1}\bar{y}_{i+1}}{\bar{y}'_{i+1}},~  y'_{i+1}=\frac{x_{i+1}y_{i}}{x'_{i}}, ~z_{i+1}=\frac{x_{i+1}\bar{y}_{i+1}z_{i}}{\bar{y}'_{i+1}x'_{i}},~{\rm if}~w_{i}\neq 0,\cr
\bar{y}'_{i+1}=0,~x'_{i}=1,~\bar{x}'_{i+1}=0, ~y'_{i+1}=x_{i+1}y_{i},~z_{i+1}=x_{i+1}\bar{y}_{i+1},~{\rm if}~w_{i}=0~{\rm and}~y_{i}\neq 0,\cr
 \bar{y}'_{i+1}=1,~x'_{i}=0,~\bar{x}'_{i+1}=\bar{x}_{i+1}\bar{y}_{i+1},~ y'_{i+1}=0,~z_{i+1}=x_{i+1}\bar{y}_{i+1},~{\rm otherwise},
 \end{cases}
\end{cases}\]
which costs at most $8m$ subtraction-free arithmetic operations.

\end{itemize}

Equipping with the preliminaries above, we are ready to  extract the representation of an arbitrary submatrix $A(\alpha|\beta)$ from the original   $ \mathcal{BD}(A)$.
We first consider the submatrix $A(r|:]$ ($A[:|r)$) by deleting the $r$th row (column) of $A\in\Rnm$. Denote
\eq{hh2.1}\mathcal{U}_{r}={\bf biupp}( \{0,1\}_{i=r}^{n})\in\Rnn,~~\mathcal{L}_{r}={\bf bilow}( \{0,1\}_{i=r}^{m})\in\Rmm.\en Then it is interesting to observe that
\[\mathcal{U}_{r}A=\mat{cc}A(r|:]\\0\rix,~~A\mathcal{L}_{r}=\mat{cc}A[:|r)&0\rix.\] Further, deleting the last row (column) of $\mathcal{U}_{r}A$ ($A\mathcal{L}_{r}$) is equivalent to left-multiplying $I_{n-1,n}$ (right-multiplying $I_{m,m-1}$), i.e.,
\[A(r|:]=I_{n-1,n}\mathcal{U}_{r}A,~~A[:|r)=A\mathcal{L}_{r}I_{m,m-1}.\] Thus,   when deleting the $r$th row (or column) of $A$, we compute the representation of the resulting submatrix $A(r|:]$ (or $A[:|r)$)   as follows:
\begin{itemize}\item first, compute the   representation of $\bar{A}=\mathcal{U}_{r}  A$ (or $\bar{A}=A\mathcal{L}_{r}$) by  using   (\ref{pt1}), (\ref{pt2}) and (\ref{pt3}),  costing at most $O((n-r)m)$ (or $O(n(m-r))$)  subtraction-free operations; \item  then, compute the  representation of   $A'=I_{n-1,n}\bar{A}$ (or $A'=\bar{A}I_{m,m-1}$)  by using   (\ref{dr}), costing at most $O(m)$  (or $O(n)$) subtraction-free operations. \end{itemize}

Now, given any submatrix $A(\alpha|\beta)$ of $A\in\Rnm$ with
$\alpha=(\alpha_{i})\in Q_{t,n}$ and $\beta=(\beta_{i})\in Q_{l,m}$,  using the argument above, it can be depicted as
\[A(\alpha|\beta)=\prod_{i=1}^{t}I_{n-t+i-1,n-t+i}\mathcal{U}_{\alpha_{i}}\cdot A\cdot \prod_{l}^{i=1}\mathcal{L}_{\beta_{i}}I_{m-l+i,m-l+i-1}.\]
Therefore, we summarize Algorithm \ref{a1} to  compute the representation of any submatrix $A(\alpha|\beta)$ in $t+l$   extracting steps, each step corresponding to the deletion of one row indexed in $\alpha$ or one column indexed in $\beta.$ The cost of the whole procedure is at most  $O((t+l)nm).$  By now, the issue of computing representations for arbitrary submatrices in general TN matrices has been addressed affirmatively, and we   enhance the existing capabilities and broaden the scope of extracting submatrix representations in nonnegative bidiagonal products, including general TN matrices.  
\begin{algorithm} \caption{Given  the representation $ \mathcal{BD}(A)$  of $A\in\Rnm$, this algorithm computes the   representation of an arbitrary submatrix $A(\alpha|\beta)$. }\label{a1}
{\bf Input:}   the representation $\mathcal{BD}(A)$,  the index sets $\alpha=(\alpha_{i})\in Q_{t,n}$ and $\beta=(\beta_{i})\in Q_{l,m}$.\\
{\bf Output:} the representation $\mathcal{BD}(A(\alpha|\beta))$.
\begin{algorithmic}[1]
\For{$i=t:-1:1$}
\State   Set $A:=\mathcal{U}_{\alpha_{i}}A $, where $\mathcal{U}_{\alpha_{i}}$ is   as in (\ref{hh2.1}).
\State Compute $\mathcal{BD}(A)$ by the procedures (\ref{pt1}), (\ref{pt2}) and (\ref{pt3}).
\State  Compute $\mathcal{BD}(A)$ of $A:=I_{n-1,n}A$ by  the downdating procedure (\ref{dr}).
\State Set $n:=n-1$.
\EndFor
\For{$i=l:-1:1$}
\State   Set $A:=A\mathcal{L}_{\beta_{i}}$, where $\mathcal{L}_{\beta_{i}}$ is   as in (\ref{hh2.1}).
\State Compute $\mathcal{BD}(A)$ by the procedures (\ref{pt1}), (\ref{pt2}) and (\ref{pt3}).
\State  Compute $\mathcal{BD}(A)$ of $A:=AI_{m,m-1}$ by  the downdating procedure (\ref{dr}).
\State Set $m:=m-1$.
\EndFor
\end{algorithmic}
\end{algorithm}


Finally, we conclude this section with how to compute  the representation $\mathcal{BD}(A)$ of a   product $A=B_1B_2\ldots B_K\in \mathbb{R}^{n_{0}\times n_{K}}$ as in (\ref{eq:biprod}). Suppose each $B_{i}\in \mathbb{R}^{n_{i-1}\times n_{i}}$ ($1\leq i\leq K$)  has $s_{i}$ nontrivial element pairs.
Then the total number of  nontrivial element pairs   is $S=\sum_{i=1}^{K}s_{i}$. We consider the case $n_{0}\geq n_{K}$ (the case  $n_{0}< n_{K}$ can be treated alike). Observe that the lower (upper) bidiagonal matrix $B_{i}$   can be rewritten as \eq{lb}B_{i}=B'_{i}I_{n_{i-1}\times n_{i}}\quad (B_{i}=I_{n_{i-1}\times n_{i}}B'_{i}),\en where $B'_{i}$ is square lower (or upper) bidiagonal having the same   nontrivial element pairs as those of $B_{i}$.
Set $A_{K}=I_{n_{K}}$. Then we  derive the  representation of $A$    by sequentially computing the    representation of  $A_{i-1}=B_{i}A_{i}$ ($i=K,\ldots,1$) as follows:
\begin{itemize} \item if $B_{i}$ is lower  bidiagonal, then by (\ref{lb}),  we first perform $A_{i}:=I_{n_{i-1}\times n_{i}}A_{i}$ by (\ref{dr2}) or (\ref{dr}),   then, compute  $A_{i-1}:=B'_{i}A_{i}$ by the procedure  (\ref{pt3}); \item if $B_{i}$ is upper  bidiagonal, then by (\ref{lb}),  we first compute  $A_{i-1}:=B'_{i}A_{i}$ by  the procedures (\ref{pt1}), (\ref{pt2}) and (\ref{pt3}),   then, perform $A_{i-1}:=I_{n_{i-1}\times n_{i}}A_{i-1}$ by (\ref{dr2}) or (\ref{dr});\end{itemize}
which costs $O(s_{i}n_{K})$  operations. Hence, the total cost   is  $O(\sum_{i=1}^{K}s_{i}n_{K})=O(S n_{K})$.

\section{The  deflation method for zero singular values}  In this section,  we develop a    method   to  exactly deflate all the zero singular values for  a nonnegative bidiagonal product of arbitrary rank.

It is important to note that deflating zero singular values exactly may not always be possible due to numerical precision limitations. The idea behind our exact deflation goes as follows: given  $ \mathcal{BD}(A)=(\{\bar{g}_{ij},g_{ij}\})\in\Rnm$,
 \begin{itemize}\item if $\bar{g}_{k1}=0$ for some $2\leq k\leq n$, then $A[k-1|:]=0$, which means that   zero rows of $A$ are implicitly revealed by  zero elements  in the first column of $\mathcal{BD}(A)$.  Thus, by the position of the zero element, we have a permutation matrix
\eq{p1}P={\bf bilow}(\{0,1\}_{i=k-1}^{n})+e_{k-1}e_{n}^{T}\in\Rnn,\en
 such that \[P^{T}A=\mat{cc}A(k-1|:]\\0\rix.\]
So, we can deflate  a possible zero singular value by deleting the zero row. It then remains to compute the submatrix representation $\mathcal{BD}(A(k-1|:])$  by Algorithm \ref{a1},  on which we  go on to  implicitly deflate zero singular values. Similar treatment applies to the case where $\bar{g}_{1k}=0$ for some $2\leq k\leq m,$ which implies  $A[:|k-1]=0;$
\item if all $\bar{g}_{i1}\neq 0$,
i.e., $\bar{g}_{i1}=1$  for all $2\leq i\leq n$, then no zero row of $A$ is revealed. In this case, denote $L={\bf bilow}(\{1,-g_{i1}\}_{i=2}^{n})$, and set $A'=LA$, then $ \mathcal{BD}(A')$ is obtained simply by setting the elements $g_{i1}=0$ ($2\leq i\leq n$) in $ \mathcal{BD}(A)$. Now, zero rows of $A'$ are revealed again by   zero elements $\bar{g}_{k2}=0$ in the second column of $ \mathcal{BD}(A')$.  Similar treatment applies to the case that $\bar{g}_{1j}\neq 0$ for all $2\leq j\leq m.$
    \end{itemize}


The other technique is to express the inverse of a lower bidiagonal matrix as the product of an orthogonal matrix and the inverse of an upper bidiagional matrix. More detailed, given  $L_{k}={\bf bilow}(\{\bar{x}_{i},-x_{i}\}_{i=k}^{n})\in\Rnn$ ($1\leq k\leq n$),   it is easy to find   Givens rotations $G_{i}\in\Rnn$ ($k+1\leq i\leq n$) such that
\[L_{k}G_{n}G_{n-1}\ldots G_{k+1}=U_{k}={\bf biupp}(\{\bar{y}_{i},-y_{i}\}_{i=k}^{n})\in\Rnn,\]
where
\[G_{i}[i-1,i]=\mat{cc}c_{i}&-s_{i}\\s_{i}&c_{i}\rix,~\forall~k+1\leq i\leq n,\]
and
\[\begin{cases}z_{n}=\bar{x}_{n},~
\bar{y}_{i}=\sqrt{z_{i}^{2}+x_{i-1}^{2}},~ c_{i}=z_{i}/\bar{y}_{i},~s_{i}=x_{i-1}/\bar{y}_{i},\cr y_{i-1}=s_{i}\bar{x}_{i-1},~z_{i-1}=c_{i}\bar{x}_{i-1},~
\bar{y}_{k}=z_{k},\end{cases}i=n,n-1,\ldots,k+1,\]
which costs at most $8(n-k)$ subtraction-free  operations. If all $\bar{x}_i$ and $x_i$ are nonnegative, then so are all $c_i,s_i,\bar{y}_i,y_i.$ Set $G=G_{n} \ldots G_{k+1}$. Then $G$ is orthogonal. Moreover, if $L_{k}$ is nonsingular, then $L^{-1}_{k}$ is expressed as
\eq{q1} L_{k}^{-1}=G\cdot U_{k}^{-1},\quad {\rm where}~ U_k^{-1}=\prod_{i=k}^nU_{i:i}~{\rm with}~U_{i:i}={\bf biupp}(\{1/\bar{y}_i,y_i/\bar{y}_i\}).\en

In what follows, we are   to  illustrate the   deflation of zero singular values for   a   nonnegative bidiagonal product $A=\prod_{i=1}^{K}B_{i}\in \mathbb{R}^{n_{0}\times n_{K}}$ as in (\ref{eq:biprod}) having   $S$ nontrivial element pairs in the factors $B_{i}\in \mathbb{R}^{n_{i-1}\times n_{i}}$ ($1\leq i\leq K$). Set the minimun dimension $r=\min_{0\leq i\leq K}\{n_{i}\}$. In Section \ref{sec4.1}, we   provide a detailed explanation of the deflation method for  the case $r=\min\{n_{0},n_{K}\}$, and  we then extend it to a general case in Section \ref{sec4.2}. Finally, the   overall SVD algorithm is established in Section \ref{sec4.3}, and error analysis is  conducted to illustrate the high relative accuracy of our proposed method.

 \subsection{The deflation for  the case $r=\min\{n_{0},n_{K}\}$} \label{sec4.1} This section is dealt with   deflating zero singular values  for the case $ r=\min\{n_{0},n_{K}\}$.

 Basing on the minimum dimension, we  compute the representation $\mathcal{BD}(A)=(\{\bar{g}_{ij},g_{ij}\})$ by using the operations (\ref{dr2}), (\ref{dr}), (\ref{pt1}), (\ref{pt2}) and (\ref{pt3}), which costs at most $O(Sr)$ subtraction-free operations.

 Focusing on  ${\cal BD}(A)=(\{\bar{g}_{ij},g_{ij}\})\in \mathbb{R}^{n_{0}\times n_{K}}$, we next illustrate how to carry out orthogonal transformations to   deflate zero singular values of $A$,  which is explained by the aid of an example with $n_{0}=5$ and $n_{K}=4$.

 As a first step, we are to deflate  zero singular values of $A$  revealed by  the  zero elements  $\bar{g}_{i1}=0$ in the first column of  $\mathcal{BD}(A)$. Assume that
 \[\bar{g}_{i_{1},1}=\bar{g}_{i_{2},1}=\ldots=\bar{g}_{i_{h},1}=0,\quad i_{1}<i_{2}<\cdots<i_{h}.\]
 Set $\alpha=(i_{1}-1,i_{2}-1,\ldots,i_{h}-1)$. Then one can detect   $A[\alpha|:]=0$. Thus, there have  permutation matrices $P_{i_{j}}$ ($1\leq j\leq h$) as in (\ref{p1}) such that $$P_{i_{1}}^{T}P_{i_{2}}^{T}\ldots P_{i_{h}}^{T}A=\mat{cc}A(\alpha|:]\\0\rix.$$
 So, we orthogonally deflate  possible $h$ zero singular values by deleting the $h$ zero  rows,   and the substantial operation is to compute    $\mathcal{BD}(A(\alpha|:])$ by Algorithm \ref{a1}, which costs at most $O(hn_{0}n_{K})$ subtraction-free operations. Moreover, in this way we can orthogonally deflate the associated zero singular values to arrive at the case where   $\bar{g}_{i1}=1$ for all $i$. For the simplicity of our description about the next operation,  in our example  $A\in \mathbb{R}^{5\times 4}$,  assume that
  \begin{eqnarray*}A&=:&\small{\begin{bmatrix}
  			\begin{smallmatrix}\{\bar{g}_{11},g_{11}\}&\{\bar{g}_{12},g_{12}\}&\{\bar{g}_{13},g_{13}\}&\{\bar{g}_{14},g_{14}\}\\
  				\{1,g_{21}\}&\{\bar{g}_{22},g_{22}\}&\{\bar{g}_{23},g_{23}\}&\{\bar{g}_{24},g_{24}\}\\
  				\{1,g_{31}\}&\{\bar{g}_{32},g_{32}\}&\{\bar{g}_{33},g_{33}\}&\{\bar{g}_{34},g_{34}\}\\
  				\{1,g_{41}\}&\{\bar{g}_{42},g_{42}\}&\{\bar{g}_{43},g_{43}\}&\{\bar{g}_{44},g_{44}\}\\
  				\{1,g_{51}\}&\{\bar{g}_{52},g_{52}\}&\{\bar{g}_{53},g_{53}\}&\{\bar{g}_{54},g_{54}\}
  	\end{smallmatrix}\end{bmatrix}}. \end{eqnarray*}

  Further, we are to perform an orthogonal transformation to zero the elements $g_{i1}$. Set $X={\bf bilow}(\{1,-g_{i+1,1}\}_{i=1}^{n_{0}-1})\in \mathbb{R}^{n_{0}\times n_{0}}$. Denote $A'=XA$. Then
  \begin{eqnarray*}A'
  	&=:&\small{\begin{bmatrix}
  			\begin{smallmatrix}\{\bar{g}_{11},g_{11}\}&\{\bar{g}_{12},g_{12}\}&\{\bar{g}_{13},g_{13}\}&\{\bar{g}_{14},g_{14}\}\\
  				\{1,0\}&\{\bar{g}_{22},g_{22}\}&\{\bar{g}_{23},g_{23}\}&\{\bar{g}_{24},g_{24}\}\\
  				\{1,0\}&\{\bar{g}_{32},g_{32}\}&\{\bar{g}_{33},g_{33}\}&\{\bar{g}_{34},g_{34}\}\\
  				\{1,0\}&\{\bar{g}_{42},g_{42}\}&\{\bar{g}_{43},g_{43}\}&\{\bar{g}_{44},g_{44}\}\\
  				\{1,0\}&\{\bar{g}_{52},g_{52}\}&\{\bar{g}_{53},g_{53}\}&\{\bar{g}_{54},g_{54}\}
  	\end{smallmatrix}\end{bmatrix}}, \end{eqnarray*} which is obtained simply by setting the elements $ g_{i1}=0$ ($2\leq i\leq n_{0}$) in $\mathcal{BD}(A)$.  Notice that the transformation $A'=XA$ is not orthogonal.  So, by \eqref{q1}, we have an orthogonal matrix $G\in \mathbb{R}^{n\times n}$ such that
  \[G^{T}X^{-1}=Y=\prod_{i=1}^{n_{0}}U_{i:i}\left(1/\bar{y}_{i},y_{i}/\bar{y}_{i}\right),~{\rm where}~\bar{y}_{i}> 0,~y_{i}\geq 0,~\forall~1\leq i\leq n_{0}.\]
  Thus, let $A^{(1)}=G^{T}A$, then
  \[A^{(1)}=(G^{T}X^{-1})(XA)=YA'=:\small{\begin{bmatrix}
  		 \begin{smallmatrix}\{\bar{g}^{(1)}_{11},g^{(1)}_{11}\}&\{\bar{g}^{(1)}_{12},g^{(1)}_{12}\}&\{\bar{g}^{(1)}_{13},g^{(1)}_{13}\}&\{\bar{g}^{(1)}_{14},g^{(1)}_{14}\}\\
  			\{1,0\}&\{\bar{g}^{(1)}_{22},g^{(1)}_{22}\}&\{\bar{g}^{(1)}_{23},g^{(1)}_{23}\}&\{\bar{g}^{(1)}_{24},g^{(1)}_{24}\}\\
  			\{1,0\}&\{\bar{g}^{(1)}_{32},g^{(1)}_{32}\}&\{\bar{g}^{(1)}_{33},g^{(1)}_{33}\}&\{\bar{g}^{(1)}_{34},g^{(1)}_{34}\}\\
  			\{1,0\}&\{\bar{g}^{(1)}_{42},g^{(1)}_{42}\}&\{\bar{g}^{(1)}_{43},g^{(1)}_{43}\}&\{\bar{g}^{(1)}_{44},g^{(1)}_{44}\}\\
  			\{1,0\}&\{\bar{g}^{(1)}_{52},g^{(1)}_{52}\}&\{\bar{g}^{(1)}_{53},g^{(1)}_{53}\}&\{\bar{g}^{(1)}_{54},g^{(1)}_{54}\}
  \end{smallmatrix}\end{bmatrix}},\] which is obtained by applying   the operations (\ref{pt1}), (\ref{pt2}) and (\ref{pt3}) to $YA'$, costing at most $O(n_{0}n_{K})$ subtraction-free  operations. By now,   possible zero singular values of $A^{(1)}$, and so $A$, can be  further revealed   by  the elements   in $\mathcal{BD}(A^{(1)})$. 

As  a second step, we go on to orthogonally deflate  zero singular values  of $A^{(1)}$ revealed by   the elements  $\bar{g}^{(1)}_{1j}=0$ in the first row of   $\mathcal{BD}(A^{(1)})$. Assume that \[\bar{g}^{(1)}_{1,j_{1}}=\bar{g}^{(1)}_{1,j_{2}}=\ldots=\bar{g}^{(1)}_{1,j_{l}}=0,\quad j_{1}<j_{2}<\cdots<j_{l}.\]
Set $\beta=(j_{1}-1,j_{2}-1,\ldots,j_{l}-1)$. Then one can detect   $A^{(1)}[:|\beta]=0$. Thus, there have   permutation matrices $Q_{j_{k}}$ ($1\leq k\leq l$) as in (\ref{p1}) such that \[A^{(1)}Q_{j_{l}}\ldots Q_{j_{2}}Q_{j_{1}}=\mat{cc}A^{(1)}[:|\beta)&0\rix.\]
So, we  orthogonally deflate  possible $l$ zero singular value by deleting the $l$ zero  columns,   and the substantial operation is to compute  $\mathcal{BD}(A^{(1)}[:|\beta))$    by Algorithm \ref{a1}, which costs at most $O(ln_{0}n_{K})$ subtraction-free operations. In this way, we can   deflate the associated zero singular values to arrive at the case that   $\bar{g}^{(1)}_{1j}=1$ for all $j$. For for the simplicity of our description about the next operation, in our example  $A^{(1)}\in \mathbb{R}^{5\times 4}$,   assume that
\[A^{(1)}=:\small{\begin{bmatrix}
  		 \begin{smallmatrix}\{\bar{g}^{(1)}_{11},g^{(1)}_{11}\}&\{1,g^{(1)}_{12}\}&\{1,g^{(1)}_{13}\}&\{1,g^{(1)}_{14}\}\\
  			\{1,0\}&\{\bar{g}^{(1)}_{22},g^{(1)}_{22}\}&\{\bar{g}^{(1)}_{23},g^{(1)}_{23}\}&\{\bar{g}^{(1)}_{24},g^{(1)}_{24}\}\\
  			\{1,0\}&\{\bar{g}^{(1)}_{32},g^{(1)}_{32}\}&\{\bar{g}^{(1)}_{33},g^{(1)}_{33}\}&\{\bar{g}^{(1)}_{34},g^{(1)}_{34}\}\\
  			\{1,0\}&\{\bar{g}^{(1)}_{42},g^{(1)}_{42}\}&\{\bar{g}^{(1)}_{43},g^{(1)}_{43}\}&\{\bar{g}^{(1)}_{44},g^{(1)}_{44}\}\\
  			\{1,0\}&\{\bar{g}^{(1)}_{52},g^{(1)}_{52}\}&\{\bar{g}^{(1)}_{53},g^{(1)}_{53}\}&\{\bar{g}^{(1)}_{54},g^{(1)}_{54}\}
  \end{smallmatrix}\end{bmatrix}}.\]
 Further, we are to orthogonally zero the  elements $g^{(1)}_{1j}$ in the first row of the representation. 
 Set $\mathcal{Y}={\bf biupp}(\{1,-g^{(1)}_{1,j+1}\}_{j=2}^{n_{K}-1})$, and denote $A'^{(1)}=A^{(1)}\mathcal{Y}$. Then
  \begin{eqnarray*}A'^{(1)}=:\small{\begin{bmatrix}
  			\begin{smallmatrix}\{\bar{g}^{(1)}_{11},g^{(1)}_{11}\}&\{1,g^{(1)}_{12}\}&\{1,0\}&\{1,0\}\\
  				\{1,0\}&\{\bar{g}^{(1)}_{22},g^{(1)}_{22}\}&\{\bar{g}^{(1)}_{23},g^{(1)}_{23}\}&\{\bar{g}^{(1)}_{24},g^{(1)}_{24}\}\\
  				\{1,0\}&\{\bar{g}^{(1)}_{32},g^{(1)}_{32}\}&\{\bar{g}^{(1)}_{33},g^{(1)}_{33}\}&\{\bar{g}^{(1)}_{34},g^{(1)}_{34}\}\\
  				\{1,0\}&\{\bar{g}^{(1)}_{42},g^{(1)}_{42}\}&\{\bar{g}^{(1)}_{43},g^{(1)}_{43}\}&\{\bar{g}^{(1)}_{44},g^{(1)}_{44}\}\\
  				\{1,0\}&\{\bar{g}^{(1)}_{52},g^{(1)}_{52}\}&\{\bar{g}^{(1)}_{53},g^{(1)}_{53}\}&\{\bar{g}^{(1)}_{54},g^{(1)}_{54}\}
  	\end{smallmatrix}\end{bmatrix}}, \end{eqnarray*} which is obtained simply by setting the elements $g^{(1)}_{1j}=0$ ($3\leq j\leq n_{K}$) in $\mathcal{BD}(A^{(1)})$. Further, by \eqref{q1}, we have an orthogonal matrix $V $ such that
  \[\mathcal{Y}^{-1}V=\mathcal{X}=\prod_{n_{K}}^{i=2}L_{i:i}\left(1/\bar{x}_{i},x_{i}/\bar{x}_{i}\right),~{\rm where}~\bar{x}_{i}> 0,~x_{i}\geq 0,~\forall ~2\leq i\leq n_{K}.\]
   Thus, let $A^{(2)}=A^{(1)}V$, then
  \[A^{(2)}=(A^{(1)}\mathcal{Y})(\mathcal{Y}^{-1}V)=A'^{(1)}\mathcal{X}=:\small{\begin{bmatrix}
  		\begin{smallmatrix}\{\bar{g}^{(2)}_{11},g^{(2)}_{11}\}&\{1,g^{(2)}_{12}\}&\{1,0\}&\{1,0\}\\
  			\{1,0\}&\{\bar{g}^{(2)}_{22},g^{(2)}_{22}\}&\{\bar{g}^{(2)}_{23},g^{(2)}_{23}\}&\{\bar{g}^{(2)}_{24},g^{(2)}_{24}\}\\
  			\{1,0\}&\{\bar{g}^{(2)}_{32},g^{(2)}_{32}\}&\{\bar{g}^{(2)}_{33},g^{(2)}_{33}\}&\{\bar{g}^{(2)}_{34},g^{(2)}_{34}\}\\
  			\{1,0\}&\{\bar{g}^{(2)}_{42},g^{(2)}_{42}\}&\{\bar{g}^{(2)}_{43},g^{(2)}_{43}\}&\{\bar{g}^{(2)}_{44},g^{(2)}_{44}\}\\
  			\{1,0\}&\{\bar{g}^{(2)}_{52},g^{(2)}_{52}\}&\{\bar{g}^{(2)}_{53},g^{(2)}_{53}\}&\{\bar{g}^{(2)}_{54},g^{(2)}_{54}\}
  \end{smallmatrix}\end{bmatrix}},\]  which is obtained by applying   the operations (\ref{pt1}), (\ref{pt2}) and (\ref{pt3}) to $A'^{(1)}\mathcal{X}$, costing at most $O(n_{0}n_{K})$ subtraction-free operations. Here,  $g^{(2)}_{12}$ is not eliminated, since attempting to eliminate it will disturb the previously reduced  element pairs $\{1,0\}$.

  Therefore, we have orthogonally deflated  some zero singular values of $A$ to get the resulting matrix
  \[A^{(2)}=G^{T}AV=\mat{cc}1&0\\0&A'^{(2)} \rix U_{1:1}(g^{(2)}_{11},g^{(2)}_{11}g^{(2)}_{12}),\]
 where $g^{(2)}_{11}\neq 0$, and
 \[\mathcal{BD}(A'^{(2)})=\small{\begin{bmatrix}
  		\begin{smallmatrix}
  			\{\bar{g}^{(2)}_{22},g^{(2)}_{22}\}&\{\bar{g}^{(2)}_{23},g^{(2)}_{23}\}&\{\bar{g}^{(2)}_{24},g^{(2)}_{24}\}\\
  			\{\bar{g}^{(2)}_{32},g^{(2)}_{32}\}&\{\bar{g}^{(2)}_{33},g^{(2)}_{33}\}&\{\bar{g}^{(2)}_{34},g^{(2)}_{34}\}\\
  			\{\bar{g}^{(2)}_{42},g^{(2)}_{42}\}&\{\bar{g}^{(2)}_{43},g^{(2)}_{43}\}&\{\bar{g}^{(2)}_{44},g^{(2)}_{44}\}\\
  			\{\bar{g}^{(2)}_{52},g^{(2)}_{52}\}&\{\bar{g}^{(2)}_{53},g^{(2)}_{53}\}&\{\bar{g}^{(2)}_{54},g^{(2)}_{54}\}
  \end{smallmatrix}\end{bmatrix}}.\] 
 More importantly, the remaining zero singular values of   $A$  are revealed by  $\mathcal{BD}(A'^{(2)})$.
Hence, proceeding with such  deflations on the representation  and going on,  finally all the zero singular values of  $A$ are orthogonally deflated, arriving at a full-rank bidiagonal matrix, e.g., our example $A\in \mathbb{R}^{5\times 4}$  is orthogonally deflated into the full-rank bidiagonal matrix
\[B=\prod _{4}^{i=1}U_{i:i}(g^{(4)}_{ii},g^{(4)}_{ii}g^{(4)}_{i,i+1})=:\small{\begin{bmatrix}
  		\begin{smallmatrix}\{\bar{g}^{(4)}_{11},g^{(4)}_{11}\}&\{1,g^{(4)}_{12}\}&\{1,0\}&\{1,0\}\\
  			\{1,0\}&\{\bar{g}^{(4)}_{22},g^{(2)}_{22}\}&\{1,g^{(2)}_{23}\}&\{1,0\}\\
  			\{1,0\}&\{1,0\}&\{\bar{g}^{(4)}_{33},g^{(4)}_{33}\}&\{1,g^{(4)}_{34}\}\\
  			\{1,0\}&\{1,0\}&\{1,0\}&\{\bar{g}^{(4)}_{44},g^{(4)}_{44}\}\\
  			\{1,0\}&\{1,0\}&\{1,0\}&\{1,0\}
  \end{smallmatrix}\end{bmatrix}}.\]
By now the deflation has been completed.

The   whole deflation process costs at most $O(\max\{ n_{0},n_{K}\}\cdot n_{0}n_{K})$ subtraction-free operations. Plus the cost of computing $\mathcal{BD}(A)$, the total cost   is at most  \[O(rS+\max\{ n^{2}_{0}n_{K},n^{2}_{K}n_{0}\}).\]

\subsection{The periodic deflation for a general case}\label{sec4.2}
Given a general  product $A\in \mathbb{R}^{n_{0}\times n_{K}}$  as in (\ref{eq:biprod}) having $S$ nontrivial element pairs, the direct method for deflating zero singular values of $A$ is that we first compute $\mathcal{BD}(A)$ and then apply the the  deflation method presented in Section \ref{sec4.1} to $\mathcal{BD}(A)$. The total cost of  the direct method is  at most  \eq{cost1}O(\min\{n_{0},n_{K}\}\cdot S+\max\{n^{2}_{0}n_{K},n^{2}_{K}n_{0}\}),\en
which is pessimistic if $n_{0}$ and $n_{K}$ are very large.  In this section, basing on the minimum dimension $r=\min_{0\leq i\leq K}\{n_{i}\}$,  we develop a periodic deflation method   such that the deflation is effectively carried out  costing at most  \[O(  r S+\max\{n^{2}_{0}r,n^{2}_{K}r\}), \] which, compared with (\ref{cost1}), is very preferable when $r\ll \min\{n_{0},n_{K}\}$.



Recalling the factor $B_{i},$ for $1\leq i\leq  K,$ of the product $A$ as in (\ref{eq:biprod}) has the dimension $n_{i-1}\times n_{i}$, we set the minimum dimension $r=n_T$, and split $A$ into $A=A_2A_1$, where
 \[A_{2}=\prod_{i=1}^{T}B_{i}\in\mathbb{R}^{n_{0}\times r},~A_{1}=\prod_{i=T+1}^{K}B_{i}\in\mathbb{R}^{r\times n_{K}}.\]
with the convention $\prod_{i}^{j}=I_{n}$ if $j<i$. Let $S_{1}$ and $S_{2}$ be the numbers of nontrivial element pairs in $A_{1}$ and $A_{2}$, respectively. Then  the total number is $S=S_{1}+S_{2}$. Instead of   $\mathcal{BD}(A)$,  we compute
$\mathcal{BD}(A_{1})$ and $\mathcal{BD}(A_{2})$, which cost at most $O(S_{1} r)$ and $O(S_{2} r)$ subtraction-free operations, respectively. So, the total cost   is  $O((S_{1}+S_{2})\cdot r)=O(S r)$.

In what follows, we   illustrate how to orthogonally deflate zero singular values of $A=A_{2}A_{1}$ focusing on
$\mathcal{BD}(A_{1})$ and $\mathcal{BD}(A_{2}).$


\begin{itemize}
\item[Step 1.] First, we delete    zero rows, indexed by $\alpha\in Q_{t_{1},r}$,  of $A_{1}=:(\{\bar{a}_{ij},a_{ij} \})$  revealed by  the zero elements $\bar{a}_{i1}=0$, and remove the corresponding columns of $A_{2}$, arriving at
     \[A_{2}A_{1}=A_{2}[:|\alpha)\cdot A_{1}(\alpha|:], \] where   $\mathcal{BD}(A_{1}(\alpha|:])$ and $\mathcal{BD}(A_{2}[:|\alpha))$ are computed by  Algorithm \ref{a1} costing at most $O(t_{1}rn_{K})$ and $O(t_{1}rn_{0})$ subtraction-free operations, respectively. In this way  we  arrive at the submatrices $A'_{1}\in\mathbb{R}^{r'\times n_{K}}$ and $A'_{2}\in\mathbb{R}^{n_{0}\times r'}$ such that
  $A=A'_{2}A'_{1}$, where $\mathcal{BD}(A'_{1})=( \{\bar{a}^{(1)}_{ij},a^{(1)}_{ij} \})$ satisfies that all $\bar{a}^{(1)}_{i1}=1$. We go on to zero the   elements $a^{(1)}_{i1}$.
  Set
$L_{1}={\bf bilow}(\{1,-a^{(1)}_{i+1,1}\}_{i=1}^{r'-1})$. Then
  \eq{def1}A= A'_{2}A'_{1}= (A'_{2}L^{-1}_{1})(L_{1}A'_{1})= A''_{2}A^{(1)}_{1},\en where \[A^{(1)}_{1}=L_{1}A'_{1}=:\mat{ccccc}\{\bar{a}^{(1)}_{11},a^{(1)}_{11} \}&\{\bar{a}^{(1)}_{12},a^{(1)}_{12} \}& \cdots&\{\bar{a}^{(1)}_{1n_{K}},a^{(1)}_{1n_{K}} \}\\\{1,0_{1}\}&\{\bar{a}^{(1)}_{22},a^{(1)}_{22} \}&\cdots&\{\bar{a}^{(1)}_{2n_{K}},a^{(1)}_{2n_{K}} \}\\ \vdots&\vdots&\vdots&\vdots\\\{1,0_{1}\}&\{\bar{a}^{(1)}_{r'2},a^{(1)}_{r'2} \}& \cdots&\{\bar{a}^{(1)}_{r'n_{K}},a^{(1)}_{r'n_{K}} \} \rix\in\mathbb{R}^{r'\times n_{K}}\]
 is  obtained simply by setting   $a^{(1)}_{i1}=0_{1}$ ($2\leq i\leq r'$) in $\mathcal{BD}(A'_{1})$,  and the   representation of   $A''_{2}=A'_{2}L^{-1}_{1}=A'_{2}\cdot\prod_{r'-1}^{i=1}L_{i:i}(1,a^{(1)}_{i+1,1})\in\mathbb{R}^{n_{0}\times r'}$ is computed by  using (\ref{pt1}), (\ref{pt2}) and (\ref{pt3}), costing at most $O(r' n_{0})$ subtraction-free operations. 
\item[Step 2.] Consider that deleting  zero rows of the resulting matrix $A''_{2}=:(\{\bar{b}_{ij},b_{ij} \})$ revealed by  the zero elements $\bar{b}_{i1}=0$,   whose number is assumed to be $t_{2}=n_{0}-n'_{0}$, is meant deflating   possible  zero singular values of  $A= A''_{2}A^{(1)}_{1}$.  Thus, by  the deflation in Section \ref{sec4.1}, we construct an orthogonal matrix $G_{1}$ to deflate  the  zero singular values,  which costs at most $O(t_{2}n_{0}r'+n'_{0}r')$ subtraction-free operations, arriving at
    \eq{def2}A^{(1)}=
(G^{T}_{1}A''_{2})A^{(1)}_{1}=A^{(1)}_{2}A^{(1)}_{1}, \en
where
\[A^{(1)}_{2}=G^{T}_{1}A''_{2}=:\mat{ccccc}\{\bar{b}^{(1)}_{11},b^{(1)}_{11} \}&\{\bar{b}^{(1)}_{12},b^{(1)}_{12} \}& \cdots&\{\bar{b}^{(1)}_{1r'},b^{(1)}_{1r'} \}\\\{1,0_{1}\}&\{\bar{b}^{(1)}_{22},b^{(1)}_{22} \}&\cdots&\{\bar{b}^{(1)}_{2r'},b^{(1)}_{2r'} \}\\ \vdots&\vdots&\vdots&\vdots\\\{1,0_{1}\}&\{\bar{b}^{(1)}_{n'_{0},2},b^{(1)}_{n'_{0},2} \}& \cdots&\{\bar{b}^{(1)}_{n'_{0},r'},b^{(1)}_{n'_{0},r'} \} \rix  \in \mathbb{R}^{n'_{0}\times r'}.\]
\item[Step 3.] Further, we delete    zero columns, indexed by $\beta\in Q_{t_{3},r'}$, of   $A^{(1)}_{2}=:(\{\bar{b}^{(1)}_{ij},b^{(1)}_{ij}\})$  revealed by  the zero elements $\bar{b}^{(1)}_{1j}=0$, and remove the corresponding rows of $A^{(1)}_{1}$, arriving at \[A^{(1)}_{2}A^{(1)}_{1}=A^{(1)}_{2}[:|\beta)\cdot A^{(1)}_{1}(\beta|:], \]where     $\mathcal{BD}(A^{(1)}_{2}[:|\beta))$ and $\mathcal{BD}(A^{(1)}_{1}(\beta|:])$  are computed by  Algorithm \ref{a1} costing at most $O(t_{3}r'n'_{0})$ and $O(t_{3}r'n_{K})$ subtraction-free operations, respectively.
    In this way we arrive at the submatrices $A'^{(1)}_{2}\in\mathbb{R}^{n'_{0}\times r''}$ and $A'^{(1)}_{1}\in\mathbb{R}^{r''\times n_{K}}$ such that $A^{(1)}=A'^{(1)}_{2}A'^{(1)}_{1}$, where
$\mathcal{BD}(A'^{(1)}_{2})=(\{\bar{b}^{(2)}_{ij},b^{(2)}_{ij}\})$ satisfies that all $\bar{b}^{(2)}_{1j}=1$. We go on to zero the  elements $b^{(2)}_{1j}$.  Set
$U_{1}={\bf biupp}(\{1,-b^{(2)}_{1,j+1}\}_{j=1}^{r''-1})$. Then
\[ A^{(1)}=(A'^{(1)}_{2}U_{1}) (U^{-1}_{1}A'^{(1)}_{1})=A^{(2)}_{2}A''^{(1)}_{1},\]
where \[A^{(2)}_{2}=A'^{(1)}_{2}U_{1}=:\mat{ccccc}\{\bar{b}^{(2)}_{11},b^{(2)}_{11} \}&\{1,0_{2} \}& \cdots&\{1,0_{2} \}\\\{1,0_{1}\}&\{\bar{b}^{(2)}_{22},b^{(2)}_{22} \}&\cdots&\{\bar{b}^{(2)}_{2r''},b^{(2)}_{2r''} \}\\ \vdots&\vdots&\vdots&\vdots\\\{1,0_{1}\}&\{\bar{b}^{(2)}_{n'_{0},2},b^{(2)}_{n'_{0},2} \}& \cdots&\{\bar{b}^{(2)}_{n'_{0},r''},b^{(2)}_{n'_{0},r''} \} \rix  \in \mathbb{R}^{n'_{0}\times r''}\]
is obtained simply by setting the elements $b^{(2)}_{1j}=0_{2}$  ($2\leq j\leq r''$) in $\mathcal{BD}(A'^{(1)}_{2})$,  and
\[A''^{(1)}_{1}=U^{-1}_{1}A'^{(1)}_{1}=:\mat{ccccc}\{\bar{a}''_{11},a''_{11} \}&\{\bar{a}''_{12},a''_{12} \}& \cdots&\{\bar{a}''_{1n_{K}},a''_{1n_{K}} \}\\\{1,0_{1}\}&\{\bar{a}''_{22},a''_{22} \}&\cdots&\{\bar{a}''_{2n_{K}},a''_{2n_{K}} \}\\ \vdots&\vdots&\vdots&\vdots\\\{1,0_{1}\}&\{\bar{a}''_{r'',2},a''_{r'',2} \}& \cdots&\{\bar{a}''_{r'',n_{K}},a''_{r'',n_{K}} \} \rix\in\mathbb{R}^{r''\times n_{K}}\]
 is computed by (\ref{pt1}), (\ref{pt2}) and (\ref{pt3}), costing at most $O(r''n_{K})$ subtraction-free operations. In particular, the previously reduced element pairs  $\{1,0_{1}\}$ of $A'^{(1)}_{1}$ and $A'^{(1)}_{2}$   are not disturbed. 
\item[Step 4.] Notice that deleting   zero columns of the resulting matrix $A''^{(1)}_{1}=:(\{\bar{a}''_{ij},a''_{ij} \})$  revealed by the zero elements $\bar{a}''_{1j}=0$,   whose number is assumed to be $t_{4}=n_{K}-n'_{K}$, is meant deflating   possible  zero singular values of $A^{(1)}=A^{(2)}_{2}A''^{(1)}_{1}$. Thus,
by   the deflation in Section \ref{sec4.1}, we construct an orthogonal matrix $V_{1}$ to deflate the    zero singular values, which costs at most $O(t_{4}r'' n_{K}+r'' n'_{K})$ subtraction-free operations, arriving at
    \eq{def3}
 A^{(2)}=A^{(2)}_{2}(A''^{(1)}_{1}V_{1})=A^{(2)}_{2}A^{(2)}_{1},\en
 where
 \[A^{(2)}_{1}=A''^{(1)}_{1}V_{1}=:\mat{ccccc}\{\bar{a}^{(2)}_{11},a^{(2)}_{11} \}&\{1,a^{(2)}_{12} \}&\{1,0_{2}\}& \cdots&\{1,0_{2} \}\\\{1,0_{1}\}&\{\bar{a}^{(2)}_{22},a^{(2)}_{22} \}&\{\bar{a}^{(2)}_{23},a^{(2)}_{23} \}&\cdots&\{\bar{a}^{(2)}_{2n'_{K}},a^{(2)}_{2n'_{K}} \}\\ \vdots&\vdots&\vdots&\vdots&\vdots\\\{1,0_{1}\}&\{\bar{a}^{(2)}_{r'',2},a^{(2)}_{r'',2} \}&\{\bar{a}^{(2)}_{r'',3},a^{(2)}_{r'',3} \}& \cdots&\{\bar{a}^{(2)}_{r'',n'_{K}},a^{(2)}_{r'',n'_{K}} \} \rix.\]
 In particular, the previously reduced element pairs  $\{1,0_{1}\}$ of $A''^{(1)}_{1}$    are not disturbed.
\end{itemize}

Therefore, we have orthogonally deflated the zero singular values of $A=A_{2}A_{1}$ to arrive at the resulting matrix
\[A^{(2)}=G_{1}^{T}AV_{1}=A^{(2)}_{2}A^{(2)}_{1}=\mat{cc}1&0\\0&A'^{(2)}_{2}\cdot A'^{(2)}_{1}\rix \cdot U_{1:1}(b^{(2)}_{11}a^{(2)}_{11},b^{(2)}_{11}a^{(2)}_{11}a^{(2)}_{12}),\]
where $b^{(2)}_{11}a^{(2)}_{11}\neq 0$, and
\[\mathcal{BD}(A'^{(2)}_{2})=\mat{ccccc}\{\bar{b}^{(2)}_{22},b^{(2)}_{22} \}&\cdots&\{\bar{b}^{(2)}_{2r''},b^{(2)}_{2r''} \}\\ \vdots&\vdots&\vdots\\\{\bar{b}^{(2)}_{n'_{0},2},b^{(2)}_{n'_{0},2} \}& \cdots&\{\bar{b}^{(2)}_{n'_{0},r''},b^{(2)}_{n'_{0},r''} \} \rix,\]and \[\mathcal{BD}(A'^{(2)}_{1})=\mat{ccccc}\{\bar{a}^{(2)}_{22},a^{(2)}_{22} \} &\cdots&\{\bar{a}^{(2)}_{2n'_{K}},a^{(2)}_{2n'_{K}} \}\\ \vdots &\vdots&\vdots\\\{\bar{a}^{(2)}_{r'',2},a^{(2)}_{r'',2} \} & \cdots&\{\bar{a}^{(2)}_{r'',n'_{K}},a^{(2)}_{r'',n'_{K}} \} \rix.\]
More importantly, the remaining zero singular values of  $A$   are revealed by   $\mathcal{BD}(A'^{(2)}_{1})$ and $\mathcal{BD}(A'^{(2)}_{2})$.
Therefore, proceeding with such  deflations on the representations and going on,
  we subsequently  get the orthogonal matrices $G'$ and $V'$ such that  the zero singular values of $A'=A'^{(2)}_{2} A'^{(2)}_{1}$ are   deflated to get that
\[B'= G'^{T}A'V'=\mat{cc} \bar{B}'&0\\0&0\rix\in\mathbb{R}^{(n'_{0}-1)\times (n'_{K}-1)},  \]
where   $\bar{B}'$ is  bidiagonal of full rank.

Now, denote $G=G_{1}\cdot {\rm diag}(1,G')$, $V=V_{1}\cdot {\rm diag}(1,V')$ and $\bar{B}={\rm diag}(1,\bar{B}')\cdot U_{1:1}(b^{(2)}_{11}a^{(2)}_{11},b^{(2)}_{11}a^{(2)}_{11}a^{(2)}_{12})$. Then
 we conclude that  all the zero singular values of  $A=A_{2}A_{1}$ are orthogonally deflated to get that
\eq{biform}B=G^{T}AV=\mat{cc} \bar{B}&0\\0&0\rix\in\mathbb{R}^{n_{0}\times n_{K}}, \en
 here, $G\in\mathbb{R}^{n_{0}\times n_{0}}$ and $V\in\mathbb{R}^{n_{K}\times n_{K}}$ are orthogonal, and $ \bar{B}$ is     bidiagonal of full rank.

 Therefore, we summarize the  periodic deflation   as Algorithm \ref{a5}. The whole deflation process consists of the steps of deleting (eliminating) one zero row (column)   for $\mathcal{BD}(A_{1})\in \mathbb{R}^{r\times n_{K}}$ and $\mathcal{BD}(A_{2})\in \mathbb{R}^{n_{0}\times r}$, each of which costs at most $O(rn_{K})$ and $O(rn_{0})$ subtraction-free operations, respectively. So, the cost of the deflation is at most
 \[O(\max\{n_{0},n_{K},r\}\cdot \max\{rn_{0},rn_{K}\}). \]
 Plus the costs of computing $\mathcal{BD}(A_{1})$ and $\mathcal{BD}(A_{2})$, the total cost is  at most \[O(rS+\max\{n^{2}_{0}r,n^{2}_{K}r\}),\]
which is preferable when $r$ is small enough.
In addition, when   accumulating  the orthogonal matrices $G$ and $V$, the costs are $O(n_{0}^{3})$ and $O(n_{K}^{3})$, respectively. 

\begin{algorithm} \caption{Given the product $A=B_{1}B_{2}\ldots B_{K}\in\mathbb{R}^{n_{0}\times n_{K}}$ as in (\ref{eq:biprod}), the algorithm orthogonally deflate  zero singular values of  $A$ into a bidiagonal matrix.}\label{a5}
{\bf Input:}     the nontrivial element pairs in the factors $B_{i}\in \mathbb{R}^{n_{i-1}\times n_{i}}$ ($1\leq i\leq K$).\\
{\bf Output:} the bidiagonal matrix $B=G^{T}AV$ as in (\ref{biform})
\begin{algorithmic}[1]
\State Set  $n=n_{0}$, $m=n_{K}$, $G=I_{n}$ and $V=I_{m}$.
\State Set  $t=1$ and $r=\min_{0\leq l\leq K}\{n_{l}\}=n_{T}$.
\State Split $A$ into $A=A_{2}A_{1}$ with $A_{1}=B_{T+1}\ldots B_{K}$ and $A_{2}=B_{1} \ldots B_{T}$.
\State First, compute $\mathcal{BD}(A_{1})=(\{\bar{g}^{(1)}_{ij}, g^{(1)}_{ij}\})$ and $\mathcal{BD}(A_{2})=(\{\bar{g}^{(2)}_{ij}, g^{(2)}_{ij}\})$.
\While{$t\leq r$}
\If{$0<T<K$}
\For{$i=r:-1:t+1$}
\If{$\bar{g}^{(1)}_{it}=0$}
\State Compute $A_{1}:=A_{1}(i-1|:]$ and $A_{2}:=A_{2}[:|i-1)$ by Alg. \ref{a1}.
\State $r:=r-1$
\EndIf
\EndFor
\State Denote $\mathcal{BD}(A_{1})=(\{\bar{g}^{(1)}_{ij},g^{(1)}_{ij}\})$,  and let $L^{(1)}_{t}={\bf bilow}(\{1,-g^{(1)}_{i+1,t}\}_{i=t}^{r-1})$.
 \State Compute  $ A_{1}:=L^{(1)}_{t}A_{1}$  simply by setting $g^{(1)}_{t+1,t}=\ldots=g^{(1)}_{rt}=0$.
 \State Compute  $A_{2}:= A_{2}(L^{(1)}_{t})^{-1}$ by (\ref{pt1}), (\ref{pt2}) and (\ref{pt3}).
\EndIf
\State Set $s=2$ if $T\neq 0$; otherwise, $s=1$.
\For{$i=n_{0}:-1:t+1$}
\If{$\bar{g}^{(s)}_{it}=0$}
\State Compute $A_{s}:=A_{s}(i-1|:]$  by Alg. \ref{a1}.
\State  $G:=G\cdot{\rm diag}(P_{i-1},I_{n-n_{0}})$, where $P_{i-1}\in \mathbb{R}^{n_{0}\times n_{0}}$ is as in (\ref{p1}).
\State $n_{0}:=n_{0}-1$.
\EndIf
\EndFor
\State Denote $\mathcal{BD}(A_{s})=(\{\bar{g}^{(s)}_{ij},g^{(s)}_{ij}\})$, and let $L^{(s)}_{t}={\bf bilow}(\{1,-g^{(s)}_{i+1,t}\}_{i=t}^{n_{0}-1})$.
\State Compute  $ A_{s}:=L^{(s)}_{t}A_{s}$  simply by setting $g^{(s)}_{t+1,t}=\ldots=g^{(s)}_{n_{0},t}=0$.
\State Compute   $(G^{(s)}_{t})^{T}(L^{(s)}_{t})^{-1}=(U^{(s)}_{t})^{-1}$ by (\ref{q1}).
\State Compute  $A_{s}:= (U^{(s)}_{t})^{-1}A_{s}$ by (\ref{pt1}), (\ref{pt2}) and (\ref{pt3}).
\State $G:=G\cdot {\rm diag}(G^{(s)}_{t},I_{n-n_{0}})$.
\If{$0<T<K$}
\For{$j=r:-1:t+1$}
\If{$\bar{g}^{(2)}_{tj}=0$}
\State Compute $A_{2}:=A_{2}[:|j-1)$ and $A_{1}:=A_{1}(j-1|:]$ by Alg. \ref{a1}.
\State $r:=r-1$
\EndIf
\EndFor
\State Denote $\mathcal{BD}(A_{2})=(\{\bar{g}^{(2)}_{ij},g^{(2)}_{ij}\})$,  and let $U^{(2)}_{t}={\bf biupp}(\{1,-g^{(2)}_{t,j+1}\}_{j=t}^{r-1})$.
 \State Compute  $ A_{2}:=A_{2}U^{(2)}_{t}$  simply by setting $g^{(2)}_{t,t+1}=\ldots=g^{(2)}_{tr}=0$.
 \State Compute  $A_{1}:= (U^{(2)}_{t})^{-1}A_{1}$ by (\ref{pt1}), (\ref{pt2}) and (\ref{pt3}).
 \EndIf
 \State Set $z=1$ if $T\neq K$, otherwise, $z=2$.
\For{$j=n_{K}:-1:t+1$}
\If{$\bar{g}^{(z)}_{tj}=0$}
\State Compute $A_{z}:=A_{z}[:|j-1)$  by Alg. \ref{a1}.
\State  $V:=V\cdot {\rm diag}(Q_{j-1},I_{m-n_{K}})$, where $Q_{j-1}\in \mathbb{R}^{n_{K}\times n_{K}}$ is as in (\ref{p1}).
\State $n_{K}:=n_{K}-1$.
\EndIf
\EndFor
\State Denote $\mathcal{BD}(A_{z})=(\{\bar{g}^{(z)}_{ij},g^{(z)}_{ij}\})$,  and let $U^{(z)}_{t}={\bf biupp}(\{1,-g^{(z)}_{t,i+1}\}_{i=t+1}^{n_{K}-1})$.
\State Compute  $ A_{z}:=A_{z}U^{(z)}_{t}$  simply by setting $g^{(z)}_{t,t+2}=\ldots=g^{(z)}_{t,n_{K}}=0$.
\State Compute   $(U^{(z)}_{t})^{-1}V^{(z)}_{t}=(L^{(z)}_{t})^{-1}$ by (\ref{q1}).
\State Compute  $A_{z}:= A_{z}(L^{(z)}_{t})^{-1}$ by (\ref{pt1}), (\ref{pt2}) and (\ref{pt3}).
\State $V:= V\cdot {\rm diag}(V^{(z)}_{t},I_{m-n_{K}})$.
\State Set  $t:=t+1$.
\EndWhile
\State The bidiagonal matrix $B=G^{T}AV={\rm diag}(A_{2}A_{1},0)\in\mathbb{R}^{n_{0}\times n_{K}}$.
\end{algorithmic}
\end{algorithm}


\subsection{The overall SVD algorithm}  \label{sec4.3}
Once all the zero singular values of $A$ are  orthogonally deflated to get  $B={\rm diag}(\bar{B},0)$ as in (\ref{biform}),   the nonzero ones are computed by applying the   LAPACK routine   {\tt DLASQ1} to  the   full-rank bidiagonal matrix $\bar{B}$ of order $\bar{r}$,  costing $O(\bar{r}^{2})$   operations  \cite{Dem90,Parlett,Parlett1}. Therefore,   we  summarize Algorithm \ref{nbp}    to    compute  the SVD of  the  product $A$ as in (\ref{eq:biprod}), which costs at most $O(rS+\max\{n^{2}_{0}r,n^{2}_{K}r\})$ arithmetic operations plus the $O(n_{K}^{3})$ and $O(n_{0}^{3})$ costs  of accumulating the left and right singular vector matrices, respectively. 

 \begin{algorithm} \caption{Given the product $A=B_{1}B_{2}\ldots B_{K}\in\mathbb{R}^{n_{0}\times n_{K}}$ as in (\ref{eq:biprod}),  the algorithm computes the SVD of $A$. }\label{nbp}
 {\bf Input:}  the nontrivial element pairs in the factors $B_{i}\in \mathbb{R}^{n_{i-1}\times n_{i}}$ ($1\leq i\leq K$).\\
 {\bf Output:} the SVD $A=U\cdot {\rm diag}(\Sigma,0)\cdot V^{T}$
\begin{algorithmic}[1]

\State  Deflate zero singular values of $A$ into $\mathcal{U}_{1}^{T}A\mathcal{V}_{1}={\rm diag}(\bar{B},0)$  by Algorithm \ref{a5}.
\State Compute the SVD $\mathcal{U}^{T}_{2}\bar{B}\mathcal{V}_{2}=\Sigma $ by the   LAPACK routine   {\tt DLASQ1}.
\State     Then the SVD $A=U\cdot {\rm diag}(\Sigma,0)\cdot V^{T}$, where  $U=\mathcal{U}_{1}\cdot{\rm diag}(\mathcal{U}_{2},I_{n_{0}-\bar{r}})\in \mathbb{R}^{n_{0}\times n_{0}}$ and $V= \mathcal{V}_{1}\cdot {\rm diag}(\mathcal{V}_{2},I_{n_{K}-\bar{r}})\in \mathbb{R}^{n_{K}\times n_{K}}$ are orthogonal.
\end{algorithmic}

\end{algorithm}

In what follows, we are   to   show that all the zero singular values are  deflated exactly, and all the nonzero ones are computed to high relative accuracy.

Remind that every single  subtraction-free operation causes at most $\mu$ ($2\mu$) relative error in any minor (singular value) of  a nonnegative bidiagonal product \cite{kov05}, where $\mu$ is the unit roundoff. In the sequel, the computed quantity wears a hat.

\begin{theorem} Let $A=B_{1}B_{2}\ldots B_{K}\in \mathbb{R}^{n_{0}\times n_{K}}$ be  as in (\ref{eq:biprod}) having $S$ nontrivial element pairs in the bidiagonal factors $B_{i}\in \mathbb{R}^{n_{i-1}\times n_{i}}$ ($1\leq i\leq K$), and let $r=\min_{0\leq k\leq K}\{n_{k}\}$. Denote by  $\sigma_{i}$ and $\hat{\sigma}_{i}$    the exact and computed nonzero descending-ordered singular values of  $A$ by Algorithm \ref{nbp}, respectively.  Set  $C= rS+\max\{n^{2}_{0}r,n^{2}_{K}r\}$.  Then all the zero singular values are  deflated exactly, and all the nonzero ones are computed as
\[  \hat{\sigma}_{i}=(1+\eta_{i})\sigma_{i},~|\eta_{i} |\leq  \frac{O(2C\mu)}{1-O(2C\mu)},\quad \forall~1\leq i\leq {\rm rank}(A).\]
\end{theorem}
\proof Spilt the product $A$ into $A=A_{2}A_{1}$ with $A_{1}\in \mathbb{R}^{r\times n_{K}}$ and $A_{2}\in \mathbb{R}^{n_{0}\times r}$. According to the procedure of Algorithm \ref{nbp},    the following statements hold.
\begin{itemize}\item First, we compute  $\mathcal{BD}(A_{1})$ and $\mathcal{BD}(A_{2})$ costing at most $O(S_{1}  r)$ and $O(S_{2}  r)$ subtraction-free operations, which cause at most $O(S_{1}  r\mu)$ and $O(S_{2}  r\mu)$ relative errors for any minor   of  the computed quantities   $ \hat{A}_{1}=:\mathcal{BD}(\hat{A}_{1}) $ and $ \hat{A}_{2}=:\mathcal{BD}(\hat{A}_{2}) $, respectively; i.e., for any $\alpha\in Q_{k,n_{0}}$, $\gamma\in Q_{k,r}$ and $\beta\in Q_{k,n_{K}}$,
\[\begin{cases}{\rm det}\hat{A}_{2}[\alpha|\gamma]=(1+\epsilon_{\alpha,\gamma}){\rm det}A_{2}[\alpha|\gamma],~|\epsilon_{\alpha,\gamma}|\leq \frac{O(S_{2}r\mu)}{1-O(S_{2}r\mu)},\cr
{\rm det}\hat{A}_{1}[\gamma|\beta]=(1+\epsilon_{\gamma,\beta}){\rm det}A_{1}[\gamma|\beta],~|\epsilon_{\gamma,\beta}|\leq \frac{O(S_{1}r\mu)}{1-O(S_{1}r\mu)},
\end{cases} \]
and there is   $O((S_{1}+S_{2})  r\mu)$ relative error for any singular value   of $\hat{A}_{2}\hat{A}_{1}$, i.e.,
\eq{er1}\sigma_{i}(\hat{A}_{2}\hat{A}_{1})=(1+\delta_{i})\sigma_{i}(A_{2}A_{1})=(1+\delta_{i})\sigma_{i},~|\delta_{i}|\leq \frac{O(2(S_{1}+S_{2})r\mu)}{1-O(2(S_{1}+S_{2})r\mu)},~\forall~i.\en
Notice that ${\rm det}A_{2}[\alpha|\gamma]\geq0 $ and ${\rm det}A_{1}[\gamma|\beta]\geq 0$. Thus,
  \begin{eqnarray*} {\rm det}(\hat{A}_{2}\hat{A}_{1})[\alpha|\beta]=(1+\epsilon_{\alpha,\beta})\cdot  {\rm det}(A_{2}A_{1})[\alpha|\beta],~|\epsilon_{\alpha,\beta}|\leq \frac{O((S_{1}+S_{2})r\mu)}{1-O((S_{1}+S_{2})r\mu)}.\end{eqnarray*}
   This means that ${\rm det}(\hat{A}_{2}\hat{A}_{1})[\alpha|\beta]=0$ if and only if ${\rm det}(A_{2}A_{1})[\alpha|\beta]=0$. So, \eq{rk1}{\rm rank}(\hat{A}_{2}\hat{A}_{1})={\rm rank}(A_{2}A_{1})={\rm rank}(A).\en
   \item We delete zero rows and the corresponding columns of $\hat{A}_{1}$ and $\hat{A}_{2}$ to get the $r'\times n_{K}$ and $n_{0}\times r'$ submatrices $A'_{1}=:\mathcal{BD}(A'_{1}) $ and $A'_{2}=:\mathcal{BD}(A'_{2}) $ by Algorithm \ref{a1},  costing at most $O((r-r')rn_{K})$ and  $O((r-r')rn_{0})$ subtraction-free operations, respectively. Then for the computed quantities   $\hat{A}'_{1}=:\mathcal{BD}(\hat{A}'_{1})$ and $\hat{A}'_{2}=:\mathcal{BD}(\hat{A}'_{2})$,
\eq{rk2}{\rm rank}(\hat{A}'_{2}\hat{A}'_{1})={\rm rank}(A'_{2}A'_{1})={\rm rank}(\hat{A}_{2}\hat{A}_{1}),\en
and \eq{er2}\sigma_{i}(\hat{A}'_{2}\hat{A}'_{1})=(1+\delta'_{i})\sigma_{i}(A'_{2}A'_{1})=(1+\delta'_{i})\sigma_{i}(\hat{A}_{2}\hat{A}_{1}),\en
where\[|\delta'_{i}|\leq \frac{O(2r(r-r')(n_{0}+n_{K})\mu)}{1-O(2r(r-r')(n_{0}+n_{K})\mu)},~\forall~i.\]
Further,  construct a bidiagonal matrix  $L_{1}$ to zero the elements in the first column of $\mathcal{BD}(\hat{A}'_{1})$,   arriving at $A^{(1)}_{1}=L_{1}\hat{A}'_{1}$  and  $A''_{2}=\hat{A}'_{2}L^{-1}_{1}$ as   in (\ref{def1}), costing at most $O( r'n_{0})$ subtraction-free operations. Then for the computed quantities $\hat{A}^{(1)}_{1}=:\mathcal{BD}(\hat{A}^{(1)}_{1})$ and $\hat{A}''_{2}=:\mathcal{BD}(\hat{A}''_{2})$,
\eq{rk3}{\rm rank}(\hat{A}''_{2}\hat{A}^{(1)}_{1})={\rm rank}(A''_{2}A^{(1)}_{1})={\rm rank}(\hat{A}'_{2}\hat{A}'_{1}),\en
and \eq{er3}\sigma_{i}(\hat{A}''_{2}\hat{A}^{(1)}_{1})=(1+\delta''_{i})\sigma_{i}(A''_{2}A^{(1)}_{1})=(1+\delta''_{i})\sigma_{i}(\hat{A}'_{2}\hat{A}'_{1}),\en
where\[|\delta''_{i}|\leq \frac{O(2r'n_{0}\mu)}{1-O(2r'n_{0}\mu)},~\forall~i.\]
\item We construct an orthogonal matrix $G_{1}$  to deflate zero singular values of $\hat{A}''_{2}\hat{A}^{(1)}_{1}$, arriving at  $A^{(1)}_{2}=G^{T}_{1}\hat{A}''_{2}$ as in (\ref{def2}), costing at most $O((n_{0}-n'_{0})r'n_{0}+r'n'_{0})$ subtraction-free operations. Then for the computed quantity $\hat{A}^{(1)}_{2}=:\mathcal{BD}(\hat{A}^{(1)}_{2})$,
    \eq{rk4}{\rm rank}(\hat{A}^{(1)}_{2}\hat{A}^{(1)}_{1})={\rm rank}(A^{(1)}_{2}\hat{A}^{(1)}_{1})={\rm rank}(\hat{A}''_{2}\hat{A}^{(1)}_{1}),\en
    and
   \eq{er4}\sigma_{i}(\hat{A}^{(1)}_{2}\hat{A}^{(1)}_{1})=(1+\delta'''_{i})\sigma_{i}(A^{(1)}_{2}\hat{A}^{(1)}_{1})=(1+\delta'''_{i})\sigma_{i}(\hat{A}''_{2}\hat{A}^{(1)}_{1}),\en
where\[|\delta'''_{i}|\leq \frac{O(2((n_{0}-n'_{0})r'n_{0}+r'n'_{0})\mu)}{1-O(2((n_{0}-n'_{0})r'n_{0}+r'n'_{0})\mu)},~\forall~i.\]
\item Now, we have obtained that by (\ref{rk2}), (\ref{rk3}) and (\ref{rk4}),
\eq{rk5}{\rm rank}(\hat{A}^{(1)}_{2}\hat{A}^{(1)}_{1})={\rm rank}(\hat{A}_{2}\hat{A}_{1}),\en
and by (\ref{er2}), (\ref{er3}) and (\ref{er4}),
\eq{er5}\sigma_{i}(\hat{A}^{(1)}_{2}\hat{A}^{(1)}_{1})=(1+\delta'''_{i})(1+\delta''_{i})(1+\delta'_{i})\sigma_{i}(\hat{A}_{2}\hat{A}_{1})=(1+\delta^{(1)}_{i})\sigma_{i}(\hat{A}_{2}\hat{A}_{1}),\en
where, let \[C_{1}=r(r-r')(n_{0}+n_{K})+r'n_{0}+(n_{0}-n'_{0})r'n_{0}+r'n'_{0},\] then
\[|\delta^{(1)}_{i}|\leq \frac{O(2C_{1}\mu)}{1-O(2C_{1}\mu)},~\forall~i.\]
    When   continuing with such deflations on  $\hat{A}^{(1)}_{1}\in \mathbb{R}^{r'\times n_{K}}$ and $\hat{A}^{(1)}_{2}\in \mathbb{R}^{n'_{0}\times r'}$ to get $A^{(2)}_{1}\in \mathbb{R}^{r''\times n'_{K}}$ and $A^{(2)}_{1}\in \mathbb{R}^{  n'_{0}\times r''}$ as in (\ref{def3}), which costs at most $O(C_{2})$ subtraction-free operations with \[C_{2}=r'(r'-r'')(n'_{0}+n_{K})+r''n_{K}+(n_{K}-n'_{K})r''n_{K}+r''n'_{K},\] using the same argument above,  we have that for
        the computed quantities  $\hat{A}^{(2)}_{1}=:\mathcal{BD}(\hat{A}^{(2)}_{1})$ and $\hat{A}^{(2)}_{2}=:\mathcal{BD}(\hat{A}^{(2)}_{2})$,
\eq{rk6}{\rm rank}(\hat{A}^{(2)}_{2}\hat{A}^{(2)}_{1})={\rm rank}(\hat{A}^{(1)}_{2}\hat{A}^{(1)}_{1}),\en
and
\eq{er6}\sigma_{i}(\hat{A}^{(2)}_{2}\hat{A}^{(2)}_{1})=(1+\delta^{(2)}_{i})\sigma_{i}(\hat{A}^{(1)}_{2}\hat{A}^{(1)}_{1}),~|\delta^{(2)}_{i}|\leq \frac{O(2C_{2}\mu)}{1-O(2C_{2}\mu)},~\forall~i.\en
\end{itemize}
Moreover, when proceeding with the deflations on $\hat{A}^{(2)}_{2}$ and $\hat{A}^{(2)}_{1}$ to  get the  bidiagonal matrix $\bar{B}$ as in (\ref{biform}), costing  at most $O(C-rS-C_{1}-C_{2})$ subtraction-free operations, we inductively get that for the computed quantity $\hat{\bar{B}}$,
\[{\rm rank}(\hat{\bar{B}})={\rm rank}(\hat{A}^{(2)}_{2}\hat{A}^{(2)}_{1}),\]
and
\[\sigma_{i}(\hat{\bar{B}})=(1+\zeta'_{i})\sigma_{i}(\hat{A}^{(2)}_{2}\hat{A}^{(2)}_{1}),~~|\zeta'_{i}|\leq \frac{O(2(C-rS-C_{1}-C_{2})\mu)}{1-O(2(C-rS-C_{1}-C_{2})\mu)},~\forall~i.\]
Finally, when applying the   LAPACK routine   {\tt DLASQ1} \cite{Parlett1} to $\hat{\bar{B}}$ of order $\bar{r}$, the computed singular values $\hat{\sigma}_{i}$ satisfy that
\[\hat{\sigma}_{i}=(1+\zeta''_{i})\sigma(\hat{\bar{B}}),~~ |\zeta''_{i}|\leq \frac{O(\bar{r}^{2}\mu)}{1-O(\bar{r}^{2}\mu)},~\forall~i.\]

Therefore, according to the analysis above, we conclude that \[{\rm rank}(\hat{\bar{B}})={\rm rank}(\hat{A}^{(2)}_{2}\hat{A}^{(2)}_{1})={\rm rank}(A)=\bar{r},\] which means that all the zero singular values are exactly deflated. Moreover, the computed $\bar{r}$ nonzero singular values satisfy that for all $1\leq i\leq \bar{r}$,
\[  \hat{\sigma}_{i}=(1+\zeta''_{i})(1+\zeta'_{i})(1+\delta^{(2)}_{i})(1+\delta^{(1)}_{i})(1+\delta_{i})\sigma_{i}=(1+\eta_{i})\sigma_{i},~ |\eta_{i}|\leq  \frac{O(2C\mu)}{1-O(2C\mu)}, \] which means that all the nonzero singular values are computed to high relative accuracy.
The result is proved. \eproof

Finally, we conclude this section  with accurately computing SVDs of arbitrary submatrices extracted from  a product of structured matrices with repeated nodes, e.g., Cauchy, Vandermonde, Cauchy-Vandermonde, Bernstein-Vandermonde matrices and so on,  irrespective  of rank deficiency and ill-conditioning.
\begin{itemize}\item First,  depict  these structured matrices  as  bidiagonal product forms, as shown in Section \ref{sec2.1};  \item Then, express arbitrary submatrices of  products of these structured matrices into   bidiagonal product forms, as shown in Section \ref{sec2.2}; \item Finally, compute the SVDs of  submatrices  by applying Algorithm \ref{nbp} to the derived bidiagonal product forms.\end{itemize}
Pleasantly, all the singular values including zero ones of arbitrary submatrices are computed to high relative accuracy, as demonstrated by   numerical experiments later.


\section{Numerical experiments} \label{sec5} In this section,   numerical experiments are performed to confirm the claimed high relative  accuracy of our proposed method. We test our method on  arbitrary submatrices extracted from  the products of  these structured   matrices with repeated nodes:
  \begin{itemize}\item the Cauchy matrix $A=(1/(x_{i}+y_{j}))\in \mathbb{R}^{ns_{1}\times ms_{2}}$ with ascending-ordered nodes
\[\begin{cases}x_{(i-1)s_{1}+1}=\ldots =x_{is_{1}}=x^{(c)}_{i},~\forall~1\leq i\leq n;\cr y_{(j-1)s_{2}+1}=\ldots =y_{js_{2}}=y^{(c)}_{j},~\forall~1\leq j\leq m;\end{cases}\]
\item the Vandermonde matrix $A=(x_{i}^{\lceil j/s_{2}\rceil -1})\in \mathbb{R}^{ns_{1}\times ms_{2}}$ with ascending-ordered nodes
\[x_{(i-1)s_{1}+1}=\ldots =x_{is_{1}}=x^{(v)}_{i},~\forall~1\leq i\leq n;\]
\item the Cauchy-Vandermonde matrix
\[A=(a_{ij})\in \mathbb{R}^{ns_{1}\times ms_{2}},~{\rm where}~a_{ij}=\begin{cases}1/(x_{i}+y_{j}),~\forall~1\leq j\leq ls_{2},\cr x_{i}^{\lceil (j-ls_{2})/s_{2}\rceil -1}~\forall~ls_{2}+1\leq j\leq ms_{2},\end{cases}\] with ascending-ordered nodes
\[\begin{cases}x_{(i-1)s_{1}+1}=\ldots =x_{is_{1}}=x^{(cv)}_{i},~\forall~1\leq i\leq n;\cr y_{(j-1)s_{2}+1}=\ldots =y_{js_{2}}=y^{(cv)}_{j},~\forall ~1\leq j\leq l;\end{cases}\]
\item the Bernstein-Vandermonde matrix
\[A=(a_{ij})\in \mathbb{R}^{ns_{1}\times ms_{2}},~{\rm where}~a_{ij}=\binom {m-1}{\lceil j/s_{2}\rceil-1}(1-x_{i})^{m-\lceil j/s_{2}\rceil}x_{i}^{\lceil j/s_{2}\rceil-1},\]
with ascending-ordered nodes
\[x_{(i-1)s_{1}+1}=\ldots =x_{is_{1}}=x^{(bv)}_{i},~\forall~1\leq i\leq n.\]\end{itemize}
Using the method described in  Section 2, we first obtain the {\it bidiagonal product} forms of the above structured matrices as well as their products, and then proceed with Algorithm \ref{nbp} to compute the SVDs.  The relative errors of the computed singular values $\hat{\sigma}_{i}$ are measured as
\[{\tt Rel.~error}(\hat{\sigma}_{i}) = |\hat{\sigma}_{i} -\sigma_{i}|/ \sigma_{i},\]
where the ``exact"  $\sigma_{i}$ are produced  by the   function {\tt SingularValueList} implemented in Mathematica with a (very expensive) 200-decimal digits arithmetic and rounded back to 16 decimal digits.

We compare the error results of our method with those of the Matlab command {\tt svd} applied to the explicit product. We need to highlight the differences in implementations between our method and the command  {\tt svd}. Our approach works on the defining nodes of these structured matrices, while the command {\tt svd} is applied to the  product explicitly computed from these structured matrices.
However, it should be noted that once the product is explicitly computed, there is no guarantee of relative accuracy for the smaller singular values  even with exact computation,  because even tiny errors in the computation of the matrix entries can lead to significant errors in the computed singular values, especially for the smaller ones.  Our main goal in comparing the relative errors   is simply to demonstrate that by carefully leveraging the structural properties, the SVDs of arbitrary submatrices of products of structured matrices can be computed with high relative accuracy, despite of ill conditioning and rank deficiency. All tests are conducted in MATLAB 7.0 using double precision arithmetic.

\begin{example} \label{exp1.1} Consider the product $A=A_{4}A_{3}A_{2}A_{1}$, where $A_{1}\in \mathbb{R}^{n_{2}s_{2}\times n_{1}s_{1}}$, $A_{2}\in \mathbb{R}^{n_{3}s_{3}\times n_{2}s_{2}}$, $A_{3}\in \mathbb{R}^{n_{4}s_{4}\times n_{3}s_{3}}$ and $A_{4}\in \mathbb{R}^{n_{5}s_{5}\times n_{4}s_{4}}$  are  Cauchy-Vandermonde, Bernstein-Vandermonde, Vandermonde and  Cauchy matrices, respectively;  whose nodes, denoted with superscript `cv', `bv', `v'  and `c', respectively,
 \[\begin{cases}x^{(cv)}_{i}=i/n_{2},~1\leq i\leq n_{2},~y^{(cv)}_{j}=j/n_{1},~1\leq j\leq l;\cr
 x^{(bv)}_{i}=1/(n_{3}-i+2)~1\leq i\leq n_{3};\cr
 x^{(v)}_{i}=1/(n_{4}-i+1),~1\leq i\leq n_{4};\cr
 x^{(c)}_{i}=1/(n_{5}-i+1),~1\leq i\leq n_{5};~y^{(c)}_{j}=(j+1)/n_{4},~1\leq j\leq n_{4}.\end{cases}\]
 Set $n_{1}=80$, $n_{2}=70$, $n_{3}=70$, $n_{4}=50$, $n_{5}=60$, $l=10$ and $s_{1}=2$, $s_{2}=3$, $s_{3}=2$, $s_{4}=4$, $s_{5}=3$. We compute  singular values of the submatrix $A'\in \mathbb{R}^{60\times 80}$ of $A$ having the $(3(i-1)+1)$th rows
and $(2(j-1)+2)$th columns for all $1\leq i\leq 60$ and $1\leq j\leq 80$.    As expected, our method  exactly identifies $10$ zero singular values, and computes all the nonzero singular values to high relative accuracy.
   In contrast, the command  {\tt svd} accurately computes only the largest singular values, failing to capture the $10$ zero singular values. The relative errors of the computed  nonzero singular values   are reported in   \tablename~\ref{ex2t1}.
\begin{table}[htbp]
\tiny
\caption{The  relative errors of  the computed singular values in Example \ref{exp1.1}.} \label{ex2t1}
\centering
\begin{tabular}{|c|c|c|c|c|}
\hline
$i$ & $\sigma_{i}(A')$&  ${\tt Rel.~error}(\hat{\sigma}_{i})$ by   Alg. \ref{nbp}  &    ${\tt Rel.~error}(\hat{\sigma}_{i})$ by {\tt svd}\\
\hline
1&1.574847649014136e+006& 8.8706e-016 &   1.1827e-015\\
 2&   2.488891921155727e+003& 3.8369e-015 & 1.7906e-014\\
  3&  1.530250850738087e-002&1.5871e-015 &   1.1686e-009\\
  4&  4.604560570301533e-007&2.7593e-015  &   1.1361e-005\\
 5&   1.230495724907674e-011&1.9694e-015   &    5.7054e+001\\
 6&   2.045994759189220e-016&4.5786e-015   &   2.5896e+006\\
 7&   2.674125770626443e-021&2.2507e-015    &   1.5249e+011\\
 8&   3.280221392098057e-026&3.3246e-015    &    1.1460e+016\\
 9&   2.792331556079233e-031&9.4095e-016    &   1.1193e+021\\
 10&   3.272383447034477e-036&8.1676e-016  &    8.6262e+025\\
 11&   2.500530956218926e-041&1.0194e-015   &   9.3835e+030\\
  12&  9.537769421320553e-046&1.3049e-015    &    4.1954e+034\\
 13&   4.707522299934931e-050& 1.0086e-015   &   8.1688e+037\\
 14&   2.903452644383030e-055&1.2476e-016     &    4.9917e+041\\
 15&   1.673616756730713e-060&4.2933e-015    &   1.4743e+046\\
 16&   1.159846521080633e-065&4.5447e-015    &   7.4582e+050\\
  17&  7.568811773397085e-071&2.7629e-015     &    5.0092e+055\\
  18&  2.547381715070572e-075&3.4688e-015   &   6.4984e+059\\
  19&  3.875607399766417e-080&4.6387e-015      &   2.1535e+064\\
  20&  2.150255912305010e-085&2.9236e-015   &     2.0886e+069\\
  21&  9.419854429953473e-091&7.1745e-015    &    2.7971e+074\\
  22&  1.002382461091948e-095&2.6549e-015   &     1.0141e+079\\
   23& 1.988962611556771e-100& 3.0624e-015    &    2.5014e+083\\
  24&  1.251313212016393e-105&1.8569e-015            &  3.2263e+088\\
  25&  7.041738853487294e-111&2.7273e-015   &   3.0746e+093 \\
  26&  7.548781771221338e-116&1.4931e-016   &  1.4836e+098\\
   27& 7.236595626968585e-121& 8.3177e-016   &   7.6894e+102\\
   28& 3.515150759856632e-126&2.4262e-015   &      1.2137e+108\\
  29&  2.224064308428203e-131&1.1252e-016     &    1.0440e+113\\
  30&  1.722231254363735e-136&2.2173e-015     &      7.3514e+117\\
   31& 9.649619609976743e-142&5.1326e-015   &      6.3227e+122\\
   32& 5.328970371601404e-147&1.8770e-015      &   1.0225e+128\\
    33&3.012241081997066e-152&4.9261e-015       &  9.0698e+132\\
   34& 1.483092890644080e-157&4.3619e-016     &     1.1362e+138\\
   35& 6.715748200773162e-163&3.3071e-015    &    1.7681e+143\\
   36& 2.762617922541912e-168&5.1113e-016     &    2.2882e+148\\
    37&9.830232713035425e-174&3.6531e-016        &    4.7691e+153\\
    38&2.950824109659177e-179& 3.4818e-016   &   9.3520e+158\\
    39&7.215403029682406e-185& 2.3538e-015     &  2.3604e+164\\
   40& 1.395795579034603e-190&5.3557e-016    &    8.9346e+169\\
   41& 2.087830631603745e-196&2.8455e-015    &   4.4184e+175\\
   42& 2.360667297832382e-202&9.2162e-015   &   1.5742e+181\\
   43& 1.970286285886377e-208&5.0459e-015    &   1.8417e+187\\
   44& 1.182321822224019e-214&5.4043e-015    &  1.7817e+193\\
   45& 4.946908994775014e-221&1.5894e-015     &     3.5330e+199\\
   46& 1.389686685088435e-227& 7.4192e-015   &    6.0298e+205\\
   47& 2.492871628238566e-234&4.8184e-015      &  3.2611e+212\\
  48&  2.652362330479972e-241&6.7796e-015     &   1.2274e+219\\
  49&  1.473958086168526e-248&8.6182e-015         &    1.6109e+226\\
   50& 3.177462705198649e-256&6.8880e-015        &   5.0078e+233\\
\hline
\end{tabular}

\end{table}
 \end{example}

\begin{example} \label{exp1.2} Consider the product $A=A_{1}A^{T}_{1}A_{1}$, where $A_{1}\in \mathbb{R}^{ns_{1}\times ms_{2}}$ is a  Cauchy-Vandermonde matrix,  whose nodes
 \[x^{(cv)}_{i}=i/2^{n-i+1},~1\leq i\leq n;~y^{(cv)}_{j}=j^{2}/2^{m-j+1},~1\leq j\leq l.\]
 Set $n=50$, $m=50$, $l=15$ and $s_{1}=3$, $s_{2}=2$. We compute the singular values of the submatrix $A'\in \mathbb{R}^{50\times 60}$ of $A$ having the $(3(i-1)+2)$th ($1\leq i\leq 50$) rows
and the columns from $21$ to $80$.    Once more, our method  computes all the nonzero singular values  to high relative accuracy, exactly returning    $20$ zero singular values.
    Conversely,  the  command  {\tt svd} is unable to accurately compute the smaller singular values, including the $20$
 zero singular values, with the correct number of significant digits.   The relative errors  of the computed   nonzero singular values   are reported in   \tablename~\ref{ex2t2}.
   \begin{table}[htbp]
\tiny
\caption{The  relative errors of  the computed singular values in Example \ref{exp1.2}.} \label{ex2t2}
\centering
\begin{tabular}{|c|c|c|c|c|}
\hline
$i$ & $\sigma_{i}(A')$&  ${\tt Rel.~error}(\hat{\sigma}_{i})$ by   Alg. \ref{nbp}  &    ${\tt Rel.~error}(\hat{\sigma}_{i})$ by {\tt svd}\\
\hline
 1&3.468887277416921e+129&1.4198e-015 &   0\\
  2&  8.018730362382941e+096&2.0820e-014 &5.9935e+004\\
  3&  4.573629685918318e+066&1.5216e-014 & 7.6397e+025\\
  4&  5.359733112254777e+040&       0 &2.5886e+038\\
   5& 3.850034487231992e+038&1.9625e-016 &1.0364e+028\\
   6& 1.219745346737479e+035&    0 & 2.4819e+026\\
   7& 1.042762368688207e+032&1.7276e-016 &8.2678e+022\\
   8& 1.930454825605820e+029&1.8226e-016 & 1.3024e+020\\
   9& 1.125637514239660e+027&3.6630e-015 &1.2180e+018\\
   10& 1.179589020167612e+018&2.8213e-015 &5.9881e+026\\
   11& 1.323455242467037e+013& 7.0837e-015 & 4.3037e+025\\
   12& 1.162171058127293e+010&2.4618e-015 &4.1844e+028\\
   13& 2.648280221329972e+006&1.2308e-015 &1.2481e+027\\
  14&  5.510081957368244e+000&1.9343e-015 &8.3542e+031\\
   15& 6.337837039281268e-002&4.3793e-016 &1.1133e+033\\
  16&  3.697851935885653e-003&1.5246e-015 &2.9211e+031\\
   17& 3.061156104412551e-006&8.5778e-015 &8.6077e+033\\
   18& 2.787053739109288e-009&6.9747e-015 &6.0813e+035\\
   19& 7.170489313678953e-012&5.2948e-015 &1.5246e+036\\
   20& 7.583241112107084e-017&3.2508e-016 &3.8804e+037\\
   21& 3.272553868178148e-020&6.6207e-015 & 7.9559e+039\\
   22& 4.510413184796948e-025&3.0541e-015 &2.4478e+043\\
   23& 5.059703390471064e-029& 1.0635e-014 &5.7334e+045\\
   24& 4.518370980853580e-035&1.1831e-015 & 6.0170e+049\\
   25& 1.018365532436303e-039&3.5242e-015 &2.5675e+053\\
   26& 1.077544585350837e-046&8.4823e-015 &3.6633e+059\\
   27& 8.872108292133314e-051&3.4784e-015 & 1.5567e+063\\
   28& 1.966010143181168e-059& 2.0242e-015 &1.9959e+070\\
  29&  2.394834026628231e-064& 1.1269e-015 &3.2862e+074\\
   30& 1.099622669639911e-075&7.8126e-016 &4.1482e+084\\
\hline
\end{tabular}

\end{table}
   \end{example}

   \begin{example} \label{exp1.3} Consider the product $A=A_{1}A^{T}_{1}A_{1}$, where $A_{1}\in \mathbb{R}^{ns_{1}\times ms_{2}}$ is a  Vandermonde matrix,  whose nodes
 \[x^{(v)}_{i}=(i+1)/(n^{2}-2i+1),~1\leq i\leq n.\]
 Set $n=50$, $m=50$ and $s_{1}=2$, $s_{2}=3$. We compute  singular values of the submatrix $A'\in \mathbb{R}^{70\times 50}$ of $A$ having the rows from $11$ to $80$ and the $(3(i-1)+2)$th ($1\leq i\leq 50$) columns.   A distinctive advantage of our method lies in its capability to exactly identify the $15$
 zero singular values, in stark contrast to the command {\tt svd} that fails to report any zeros.  The  relative errors of the computed nonzero  singular values   are reported in   \tablename~\ref{ex2t3}, highlighting the high relative accuracy of our method.
    \begin{table}[htbp]
\tiny
\caption{The  relative errors of  the computed singular values in Example \ref{exp1.3}.} \label{ex2t3}
\centering
\begin{tabular}{|c|c|c|c|c|}
\hline
$i$ & $\sigma_{i}(A')$&  ${\tt Rel.~error}(\hat{\sigma}_{i})$ by   Alg. \ref{nbp}  &    ${\tt Rel.~error}(\hat{\sigma}_{i})$ by {\tt svd}\\
\hline
 1&2.510397022449398e+003&  3.6229e-016 &  1.8115e-016\\
  2&  3.816985499740374e-004&1.4202e-016 &2.8474e-012\\
  3&  4.143858773645414e-011&1.2476e-015 & 3.4959e-006\\
   4& 4.242521375454278e-018&1.0895e-015 &8.1961e-001\\
  5&  4.239056395592293e-025&2.3831e-015 &1.9518e+005\\
   6& 4.168659664179253e-032&2.3636e-015 &3.8671e+010\\
   7& 4.043646122116740e-039& 2.7433e-015 & 2.6793e+015\\
   8& 3.869430900719566e-046&2.6134e-015 & 5.7379e+020\\
   9& 3.649829639076117e-053& 2.0325e-015 &7.9162e+025\\
   10& 3.389249221406257e-060& 5.2185e-015 & 1.5208e+031\\
   11& 3.093670547436131e-067&4.2596e-016 & 2.6090e+036\\
   12& 2.770991249357053e-074&6.9447e-015 &7.7554e+041\\
   13& 2.430919688798274e-081&5.9702e-015 &1.3220e+047\\
  14&  2.084507385475087e-088& 2.2758e-015 &3.1006e+052\\
  15&  1.743397226830317e-095&9.5404e-016 & 5.7334e+057\\
  16&  1.418889627711162e-102&5.7294e-015 &8.9102e+062\\
   17& 1.120959850392690e-109& 2.9520e-015 & 2.3207e+068\\
   18& 8.573784180593364e-117&3.7794e-015 &6.4604e+073\\
   19& 6.330803127549630e-124&1.2336e-014 &1.7316e+079\\
   20& 4.498919337915358e-131&6.5640e-015 &5.5462e+084\\
   21& 3.066611613512485e-138& 5.0588e-015 &1.2460e+090\\
   22& 1.997581004819857e-145& 4.6289e-015 &3.5820e+095\\
   23& 1.238413047998368e-152&1.3522e-014 & 1.2220e+101\\
   24& 7.273625411986728e-160&1.7371e-015 &3.6959e+106\\
  25&  4.026307600418828e-167&4.8632e-015 & 1.5156e+112\\
   26& 2.088084609674866e-174& 3.2246e-015 &3.5216e+117\\
   27& 1.007534076170427e-181&5.5766e-015 & 1.5141e+123\\
   28& 4.485972259301394e-189&5.3325e-016 &9.1815e+128\\
  29&  1.824599254309977e-196& 1.1331e-014 &4.4799e+134\\
   30& 6.694074753934959e-204&2.5391e-015 &2.6863e+140\\
  31&  2.178711168080216e-211& 3.8750e-015 &1.3847e+146\\
 32&   6.146438129885899e-219&6.3450e-015 &7.7022e+151\\
  33&  1.450862706608087e-226&4.9099e-015  &7.6456e+157\\
   34& 2.692622897706825e-234&1.5769e-014 &9.1431e+163\\
   35& 3.389889240846505e-242& 2.7014e-015 &6.3336e+169\\
\hline
\end{tabular}
\end{table}
   \end{example}


   \begin{example} \label{exp1.4} Consider the product $A=A_{1}A^{T}_{1}A_{1}$, where $A_{1}=:\mathcal{BD}(A_{1})=( \bar{g}_{ij},g_{ij})\in \mathbb{R}^{90\times 50}$, here all $g_{ij}=r_{i}/j$ with $r_{i}$ are randomly chosen by the Matlab command {\tt rand}, and all $\bar{g}_{ij}$ are randomly chosen by the Matlab command  {\tt randint}.
 The rank of $A$ is unknown. The explicit matrix $A_{1}$ is formed in a subtraction-free manner, ensuring a high level of relative accuracy. Since $\sigma_{i}(A)=(\sigma_{i}(A_{1}))^{3}$ for all $i$,  the command {\tt svd}  computes singular values of $A$  by  $({\tt svd}(A_{1}))^{3}$.
 As observed,  Alg. \ref{nbp}  computes all the $14$ nonzero singular values  to high relative accuracy,  and all the $36$ zero singular values are exactly returned.   However,     ${\tt (svd(A_{1}))^{3}}$  only identifies   $1$ zero singular values and exhibits a loss of relative accuracy for some small non-zero singular values. The relative errors  of the computed  nonzero singular values  are reported in   \tablename~\ref{ex2t4}.
   \begin{table}[htbp]
\tiny
\caption{The  relative errors of  the computed singular values in Example \ref{exp1.4}.} \label{ex2t4}
\centering
\begin{tabular}{|c|c|c|c|c|}
\hline
$i$ & $\sigma_{i}(A')$&  ${\tt Rel.~error}(\hat{\sigma}_{i})$ by   Alg. \ref{nbp}  &    ${\tt Rel.~error}(\hat{\sigma}_{i})$ by {\tt svd}\\
\hline
 1&8.332062400159658e+171& 1.0299e-015 &3.0896e-015\\
   2& 2.564548603621145e+089& 1.1023e-015 & 9.2338e+033\\
  3&  1.099405131386926e+038&3.4363e-016 &1.9317e+080\\
  4&  2.878741420489996e+006&1.6176e-016 &1.2702e+109\\
   5& 7.418482246460084e-005&1.0961e-015 &3.1706e+115\\
   6& 2.915960871453962e-014&6.4928e-016 &5.0791e+123\\
   7& 9.291825752068523e-016&4.4572e-015 & 1.1561e+121\\
   8& 3.041046489601403e-022& 1.5462e-015 &2.8172e+120\\
   9& 3.921877379377079e-035&4.0890e-015 &4.2790e+126\\
   10& 9.631567808072229e-052&3.0809e-016 &9.4277e+140\\
   11& 1.612410893623431e-055&6.8518e-015 &2.0076e+139\\
   12& 3.747225733496450e-059&1.8880e-015 & 1.7617e+131\\
   13& 6.896108910152905e-078& 4.1711e-016 &8.2239e+149\\
   14& 6.455473653011518e-121&4.5289e-015 &2.8651e+186 \\
\hline
\end{tabular}
\end{table}

   \end{example}

\begin{appendices}

\section{The proofs of the passing-through operations (\ref{pt1}) and (\ref{pt3})}\label{secA1}

{\bf The proof of the operation (\ref{pt1})}: Since both   $U_{1:n}L_{1:n}$ and $\bar{L}_{1:n}\bar{U}_{1:n}$ are tridiagonal, comparing the entries on the diagonal, subdiagonal and superdiagonal gives
\[\begin{cases} \bar{y}'_{i}\bar{x}'_{i}+x'_{i-1}y'_{i-1}=\bar{y}_{i}\bar{x}_{i}+x_{i}y_{i},\cr
	\bar{y}'_{i}x'_{i}=\bar{y}_{i+1}x_{i},~\bar{x}'_{i}y'_{i}=\bar{x}_{i+1}y_{i},\end{cases}~\forall~i=1,\ldots,n,\]
with the convention $x'_{0}=y'_{0}=0$.  To avoid subtraction, we introduce the intermediate quantities
$$
z_{i}=\bar{y}_{i}\bar{x}_{i}-x'_{i-1}y'_{i-1},~1\leq i\leq n
$$
and let
$$
w_{i}:=\bar{y}'_{i}\bar{x}'_{i}=z_{i}+x_{i}y_{i}.
$$
Observe that   $ z_{1}=\bar{y}_{1}\bar{x}_{1}\geq 0$.  Assume that  $z_{i}\geq 0$ and thus $w_i\ge 0.$
Then,
\begin{itemize}\item if  $w_{i}>0$,  then we   set
	\[\bar{x}'_{i}=1,~\bar{y}'_{i}=w_{i}>0,~x'_{i}=\frac{\bar{y}_{i+1}x_{i}}{\bar{y}'_{i}}\geq  0,~y'_{i}=\frac{y_{i}\bar{x}_{i+1}}{\bar{x}'_{i}}\geq 0,\]  thus, \[z_{i+1}=\bar{y}_{i+1}\bar{x}_{i+1}-x'_{i}y'_{i}=\frac{\bar{y}_{i+1}\bar{x}_{i+1}(\bar{y}'_{i}\bar{x}'_{i}-x_{i}y_{i})}{\bar{y}'_{i}\bar{x}'_{i}}
	=\frac{\bar{y}_{i+1}\bar{x}_{i+1}z_{i}}{\bar{x}'_{i}\bar{y}'_{i}}\geq 0;\]
	\item if $w_{i}=0$,  then $x_{i}y_{i}=0.$  If  $y_{i}=0$,   then we set
	\[\bar{x}'_{i}=w_{i}=0,~\bar{y}'_{i}= 1,~x'_{i}= \bar{y}_{i+1}x_{i} \geq 0,~y'_{i}=0,\]thus, \[z_{i+1}=\bar{y}_{i+1}\bar{x}_{i+1}-x'_{i}y'_{i}=\bar{y}_{i+1}\bar{x}_{i+1}\geq 0;\]
	and if  $y_{i}\neq 0$, then  $x_{i}=0$,    we set
	\[\bar{x}'_{i}=1,~\bar{y}'_{i}=w_{i}= 0,~x'_{i}=0,~y'_{i}=y_{i}\bar{x}_{i+1}\geq 0,\]thus,\[~z_{i+1}=\bar{y}_{i+1}\bar{x}_{i+1}-x'_{i}y'_{i}=\bar{y}_{i+1}\bar{x}_{i+1}\geq 0.\]
\end{itemize}
 Inductively we have shown that all the involved quantities are nonnegative, which costs at most $8n$ subtraction-free  operations. Thus,   (\ref{pt1}) is proved. \eproof

{\bf The proof of the operation (\ref{pt3})}: By comparing the $(i,i+1)$th, $(i+1,i+1)$th and $(i,i+2)$th entries on both sides of $U_{1:m}U'_{1:m}=\bar{U}'_{1:m}\bar{U}_{2:m}$, we have that $\bar{x}'_{1}=\bar{y}_{1}\bar{x}_{1}$, and
\[\begin{cases}
y'_{i}\bar{x}'_{i}+\bar{y}'_{i+1}x'_{i}=x_{i}\bar{y}_{i}+\bar{x}_{i+1}y_{i},\cr
\bar{x}'_{i+1}\bar{y}'_{i+1}=\bar{x}_{i+1}\bar{y}_{i+1},~y'_{i+1}x'_{i}=x_{i+1}y_{i},\end{cases}\forall~i=1,\ldots,m-1,\] with the convention $y'_{1}=0$.
   To avoid  subtraction, we introduce the intermediate quantities
 \[~z_{i}=x_{i}\bar{y}_{i}-y'_{i}\bar{x}'_{i},~1\leq i\leq m-1,\]
and let \[w_{i}:=\bar{y}'_{i+1}x'_{i}=z_{i}+\bar{x}_{i+1}y_{i}.\]
Observe that $z_{1}=x_{1}\bar{y}_{1}\geq 0$. Assume that $z_{i}\geq 0$ and thus $w_{i}\geq 0$. Then,
\begin{itemize}
 \item if $w_{i}>0$, then we set \[x'_{i}=w_{i}>0,~\bar{y}'_{i+1}=1,~\bar{x}'_{i+1}=\frac{\bar{x}_{i+1}\bar{y}_{i+1}}{\bar{y}'_{i+1}}\geq 0,~y'_{i+1}=\frac{x_{i+1}y_{i}}{x'_{i}}\geq 0,\] thus, \begin{eqnarray*}z_{i+1}&=&x_{i+1}\bar{y}_{i+1}-y'_{i+1}\bar{x}'_{i+1}=\frac{x_{i+1}\bar{y}_{i+1}(\bar{y}'_{i+1}x'_{i}-\bar{x}_{i+1}y_{i})}{\bar{y}'_{i+1}x'_{i}}=\frac{x_{i+1}\bar{y}_{i+1}z_{i}}{\bar{y}'_{i+1}x'_{i}}\geq 0;\end{eqnarray*}
     \item if $w_{i} =0$,  then $\bar{x}_{i+1}y_{i}=0$. If  $y_{i}= 0$,   then we set
       \[x'_{i}=w_{i}=0,~\bar{y}'_{i+1}=1,~\bar{x}'_{i+1}=\bar{x}_{i+1}\bar{y}_{i+1}\geq 0,~y'_{i+1}=0,\]thus, \[z_{i+1}=x_{i+1}\bar{y}_{i+1}-y'_{i+1}\bar{x}'_{i+1}=x_{i+1}\bar{y}_{i+1}\geq 0;\]
      and if $y_{i}\neq 0$, then   $\bar{x}_{i+1}=0$,   then we set
      \[x'_{i}=1,~\bar{y}'_{i+1}=w_{i}=0,~\bar{x}'_{i+1}=0,~y'_{i+1}=x_{i+1}y_{i}\geq 0,\]thus, \[z_{i+1}=x_{i+1}\bar{y}_{i+1}-y'_{i+1}\bar{x}'_{i+1}=x_{i+1}\bar{y}_{i+1}\geq 0.\]
     \end{itemize}
  Inductively we have shown that all the involved quantities are nonnegative, which costs at most $8m$ subtraction-free  operations.  Thus,  (\ref{pt3}) is proved. \eproof




\end{appendices}



{\bf Funding} The work of  Rong Huang is supported by   the National Natural Science Foundation  of China (Grant No.  12271153). The work of   Jungong Xue is supported by the National Science Foundation of China (Grant No. 11771100) and the Laboratory of Mathematics for Nonlinear Science, Fudan University.


{\bf Data Availability} The datasets generated during the current study are available from the corresponding authors
on reasonable request.

\section*{Declarations} 

{\bf Conflict of interest} The authors have no relevant financial or non-financial interests to disclose.

\end{document}